\date{21 August 2006}
\documentclass[reqno,11pt]{amsart}
\usepackage{latexsym,amsmath,amsfonts,amscd,amssymb}

\usepackage{eepic}
\theoremstyle{plain}
\newtheorem{theorem}{Theorem}[section]
\newtheorem{corollary}[theorem]{Corollary}
\newtheorem{lemma}[theorem]{Lemma}
\newtheorem{proposition}[theorem]{Proposition}
\theoremstyle{definition}
\newtheorem{definition}[theorem]{Definition}
\newtheorem{remark}[theorem]{Remark}
\numberwithin{equation}{section}
\newenvironment{dem}{\noindent{\bf Proof:}\ }{\hfill$\Box$\newline}

\newcommand{\s}{\sigma}

\renewcommand{\t}{\tau}

\newcommand{\cE}{\mathcal{E}}

\newcommand{\cN}{\mathcal{N}}

\newcommand{\cO}{\mathcal{O}}

\newcommand{\cS}{\mathcal{S}}
\newcommand{\cT}{\mathcal{T}}

\newcommand{\cV}{\mathcal{V}}

\newcommand{\QQ}{\mathbb{Q}}
\newcommand{\RR}{\mathbb{R}}
\newcommand{\CC}{\mathbb{C}}
\newcommand{\HH}{\mathbb{H}}
\newcommand{\PP}{\mathbb{P}}
\newcommand{\ZZ}{\mathbb{Z}}

\newcommand{\inc}{\hookrightarrow}
\newcommand{\lto}{\longrightarrow}
\newcommand{\surj}{\twoheadrightarrow}

\newcommand{\x}{\times}
\newcommand{\ox}{\otimes}
\newcommand{\iso}{\cong}

\newcommand{\frM}{{\mathfrak{M}}}
\newcommand{\coeff}{\mathop{\mathrm{coeff}}}
\newcommand{\Sym}{\mathrm{Sym}}
\newcommand{\GCD}{\mathrm{gcd}}
\newcommand{\Jac}{\mathrm{Jac}}
\DeclareMathOperator{\rk}{rk} 
 \DeclareMathOperator{\Hom}{Hom}
\DeclareMathOperator{\Ext}{Ext} \DeclareMathOperator{\End}{End}
\DeclareMathOperator{\Gr}{Gr}

\newcommand{\Res}{\mathop\mathrm{Res}}

\newcommand{\scp}{{\s_c^+}}
\newcommand{\smp}{{\s_m^+}}
\newcommand{\scm}{{\s_c^-}}
\newcommand{\scpm}{{\s_c^\pm}}
\newcommand{\moduli}{\cN_\s}

\newcommand{\Smas}{\cS_{\scp}}
\newcommand{\Smenos}{\cS_{\scm}}

\newcommand{\modulimmas}{\cN_{\smp}}

\title{Hodge polynomials of the moduli spaces of pairs}

 \thanks{First and third authors partially
supported through grant MCyT (Spain) MTM2004-07090-C03-01/02. }

\author{V. Mu\~noz}
\author{D. Ortega}
\author{M.J. V{\'a}zquez-Gallo}

\begin{document}

\begin{abstract}
 Let $X$ be a smooth projective curve of genus $g\geq 2$ over the
 complex numbers. A holomorphic pair on $X$ is a couple
 $(E,\phi)$, where $E$ is a holomorphic bundle over $X$ of rank $n$
 and degree $d$, and $\phi\in H^0(E)$ is a holomorphic section.
 In this paper, we determine the Hodge polynomials of the moduli
 spaces of rank $2$ pairs, using the theory of mixed Hodge
 structures. We also deal with the case in which $E$ has
 fixed determinant.
\end{abstract}

\maketitle

%%%%%%%%%%%%%%%%%%%%%%%%%%%%%%%%%%%%%%
\section{Introduction}
\label{sec:introduction}
%%%%%%%%%%%%%%%%%%%%%%%%%%%%%%%%%%%%%

Let $X$ be a smooth projective curve of genus $g\geq 2$ over the
field of complex numbers. Let $M(n,d)$ be the moduli space of
polystable vector bundles over $X$ of rank $n$ and degree $d$.
Similarly, we denote by $M(n,\Lambda)$ the moduli space of
polystable vector bundles of rank $n$ and determinant isomorphic
to a fixed line bundle $\Lambda$ of degree $d$ on $X$. In
\cite{DR} Desale and Ramanan gave an inductive method to determine
the Betti numbers of $M(n,\Lambda)$ for $n$ and $d$ coprime.
Later, Earl and Kirwan \cite{EK} extended the method to determine
the Hodge numbers of $M(n,d)$ and $M(n,\Lambda)$, for $n$ and $d$
coprime.

A \textit{holomorphic pair} on $X$ consists of a couple
$(E,\phi)$, where $E$ is a holomorphic bundle over $X$ of rank $n$
and degree $d$, and $\phi\in H^0(E)$ is a non-zero holomorphic
section. There is a concept of stability for a pair which depends
on the choice of a parameter $\t \in \RR$. Denote by
$\frM_{\t}(n,d)$ the moduli space of $\t$-polystable holomorphic
pairs, and $\frM_{\t}(n,\Lambda)$ the moduli space of
$\t$-polystable holomorphic pairs $(E,\phi)$, where the
determinant of $E$ is isomorphic to a fixed line bundle $\Lambda$
of degree $d$ on $X$. These moduli spaces, studied in \cite{B, BD,
GP}, are smooth for any rank $n$ at $\tau$-stable points.

In \cite{Th} Thaddeus studied how the moduli spaces
$\frM_\t(2,\Lambda)$ change when varying the parameter $\t$. The
range of the parameter $\t$ is an interval $J\subset \RR$ split by
a finite number of \textit{critical values} $\t_c$ in such a way
that, when $\t$ moves without crossing a critical value,
$\frM_{\t}(2,\Lambda)$ remains unchanged. When $\t$ crosses a
critical value, $\frM_{\t}(2,\Lambda)$ undergoes a flip (in the
sense of Mori theory), consisting of blowing-up an embedded
subvariety and then blowing-down the exceptional divisor in a
different way. The study of this process is what allows us to
obtain information on the topology of all moduli spaces
$\frM_{\t}(2,\Lambda)$, for any $\t$, once one knows such
information for one particular $\frM_{\t}(2,\Lambda)$ (usually the
one corresponding to the minimum or maximum possible value of the
parameter). This was used in \cite{Th} to compute the Poincar\'e
polynomial (i.e.,\ the Betti numbers) of $\frM_\t(2,\Lambda)$. The
argument requires an explicit construction of the blow-up and
blow-down centers of the flips, as well as the corresponding
normal bundles.

In this paper, we determine the Hodge polynomials (i.e.,\ the
Hodge numbers) of the moduli spaces $\frak M_{\t}(2,d)$ and $\frak
M_{\t}(2,\Lambda)$ without describing these blow-ups and
blow-downs. We use the theory of mixed Hodge structures introduced
by Deligne \cite{De}.

In order to set up the theory in a more general framework,  we
work with triples instead of pairs. A \textit{holomorphic triple}
$T = (E_{1},E_{2},\phi)$ on $X$ consists of two holomorphic vector
bundles $E_{1}$ and $E_{2}$ over $X$ (of ranks $n_1$ and $n_2$,
and degrees $d_1$ and $d_2$, respectively) and a holomorphic map
$\phi \colon E_{2} \to E_{1}$. There is a concept of stability for
a triple which depends, as for pairs, on the choice of a parameter
$\s \in \RR$. This gives a collection of moduli spaces
$\moduli(n_1,n_2,d_1,d_2)$, which have been studied in
\cite{BGP,BGPG,GPGM}. Again, the range of the parameter $\s$ is an
interval $I\subset \RR$ split by a finite number of critical
values. The moduli spaces change when we cross a critical value by
the removal of a subvariety and the insertion of a new subvariety.
By analogy with the situation of pairs, we call this process a
\textit{flip}, although it does not consist of a blow-up and a
blow-down in general. The subvarieties that are removed and
inserted are called the \textit{flip loci}.

When the rank of $E_2$ is one, we recover the moduli spaces of
pairs. More precisely, there is an isomorphism
 $$
 \moduli(2,1,d_1,d_2)\cong \frM_\tau(2, d_1-2d_2)\times \Jac^{d_2}X,
 $$
where $\tau$ and $\sigma$ are related by $\tau= \frac13 (\s+d_1-2
d_2)$.

Our main result is the following:

\begin{theorem} \label{thm:main}
Let $X$ be a smooth projective curve of genus $g\geq 2$. Let
$\t\in J$ be a non-critical value for the moduli space
$\frM_{\t}(2,d)$. Then $\frM_{\t}(2,d)$ is a smooth projective
variety with Hodge polynomial
 $$
  e(\frM_{\t}(2,d))= \coeff_{x^0}
  \left[\frac{(1+u)^{g}(1+v)^{g}(1+ux)^{g}(1+vx)^{g}}{(1-uv)(1-x)(1-uvx)x^{d-1-[\t]}}
    \Bigg(\frac{(uv)^{d-1-[\t]}}{1-(uv)^{-1}x}-
    \frac{(uv)^{g+1-d+2[\t]}}{1-(uv)^2x}\Bigg)\right] .
 $$
For a fixed line bundle $\Lambda$ on $X$ of degree $d$ and
non-critical value $\t\in J$, the moduli space
$\frM_{\t}(2,\Lambda)$ is a smooth projective variety with Hodge
polynomial
 $$
  e(\frM_{\t}(2,\Lambda))= \coeff_{x^0}
  \left[\frac{(1+ux)^{g}(1+vx)^{g}}{(1-uv)(1-x)(1-uvx)x^{d-1-[\t]}}
    \Bigg(\frac{(uv)^{d-1-[\t]}}{1-(uv)^{-1}x}-
    \frac{(uv)^{g+1-d+2[\t]}}{1-(uv)^2x}\Bigg)\right] .
 $$
\end{theorem}

\medskip

In Section \ref{sec:hodgepolynomials} we review the basics of
Hodge-Deligne theory. Sections~\ref{sec:stable-triples} and
\ref{sec:extensions-of-triples} recall standard results on
$\s$-stable triples used throughout the paper and give the
description of the flip loci. We particularize to triples of rank
$(n_1,n_2)=(2,1)$ in Section \ref{sec:triples(2,1)}. In Section
\ref{sec:Hodge(2,1)} we perform the computation of the Hodge
polynomials of the moduli spaces of triples of ranks $(2,1)$ and
give a proof of Theorem \ref{thm:main}. For completeness, Section
\ref{sec:Hodge(1,2)} deals with triples of rank $(1,2)$. Finally,
in Section \ref{sec:poly-bundles-rk-two}, we use the Hodge
polynomial of the moduli spaces of triples to recover the (already
known) Hodge polynomial of the moduli space of rank $2$ stable
vector bundles of odd degree.

\textbf{Acknowledgements:} First author would like to thank Peter
Gothen and Marina Logares for discussions. We are grateful to the
referee for a careful reading of the manuscript.

%%%%%%%%%%%%%%%%%%%%%%%%%%%%%%%%%%%%%%%%%%%%%%%%%%%%%%%%%%%%%%%%%%%%%%%%%%%%
\section{Hodge Polynomials}
\label{sec:hodgepolynomials}
%%%%%%%%%%%%%%%%%%%%%%%%%%%%%%%%%%%%%%%%%%%%%%%%%%%%%%%%%%%%%%%%%%%%%%%%%%%

Let us start by recalling the Hodge-Deligne theory of algebraic
varieties over $\CC$. Let $H$ be a finite-dimensional complex
vector space. A {\em pure Hodge structure of weight $k$} on $H$
%consists of a descending filtration $F$ on $H$, such that
is a decomposition
 $$
 H=\bigoplus\limits_{p+q=k} H^{p,q}
 $$
such that $H^{q,p}=\overline{H^{p,q}}$, the bar denoting complex
conjugation in $H$. (A Hodge structure is usually defined over the
field of the rational numbers, meaning that we have a vector space
$H$ over $\QQ$ whose complexification $H_\CC=H\otimes \CC$ is
endowed with a Hodge decomposition as above. Since we shall not
need this, we limit ourselves to using Hodge structures over
$\CC$.) The $k$-dimensional cohomology $H^k(Z)$ of a compact
K{\"a}hler manifold $Z$ (in particular an algebraic smooth projective
variety) has a natural Hodge structure of weight $k$ (here we
denote $H^k(Z)$ for the cohomology with complex coefficients,
although the result also holds for the cohomology with rational
coefficients). We denote
  $$
  h^{p,q}(Z)=\dim H^{p,q}\ ,
  $$
which is called the Hogde number of type $(p,q)$. A Hodge structure of
weight $k$ on $H$ gives rise to the so-called {\em Hodge filtration} $F$
on $H$, where
  $$
  F^p= \bigoplus\limits_{s\geq p} H^{s,k-s}\ ,
  $$
which is a descending filtration. Note that $\Gr_F^p H = F^p/F^{p+1}=
H^{p,q}$.

Let $H$ be a complex vector space. A {\em (mixed) Hodge structure} over
$H$ consists of an ascending weight filtration $W$ on $H$ and a descending
Hodge filtration $F$ on $H$ such that $F$ induces a pure Hodge filtration
of weight $k$ on each $\Gr^W_k H= W_k/W_{k-1}$. Again we define
  $$
  h^{p,q}(H)=\dim \, H^{p,q}\, , \qquad \text{where}\quad H^{p,q}
  =\Gr_F^p\Gr^W_{p+q} H\, .
  $$

A morphism between two Hodge structures $H$, $H'$ is a linear map
$f:H\to H'$ compatible with both the weight and Hodge filtrations.
The most relevant fact is that morphisms of Hodge structures are
strictly compatible with the filtrations $W$ and $F$ (see
\cite{De}). From this it follows that the functor $H\mapsto
H^{p,q}$ is exact.

Deligne has shown \cite{De} that, for each complex algebraic
variety $Z$, the cohomology $H^k(Z)$ carries a natural Hodge
structure, with weight filtration $W$ and Hodge filtration $F$,
which coincides with the classical (pure) Hodge structure when $Z$
is a smooth projective variety. The associated graded objects of
these filtrations are denoted $\Gr^W$ and $\Gr_F$, respectively.
We define the following Euler characteristics, by using these
filtrations:
  $$
  \chi_{p,q}(Z) = \sum_k (-1)^k \dim \Gr^p_F \Gr_{p+q}^W H^k(Z).
  $$

The cohomology groups with compact support $H_c^k(Z)$ also carry a
natural mixed Hodge structure \cite{De}, which allows one to
define $\chi^c_{p,q}$, using $H_c^k(Z)$ instead of $H^k(Z)$. If
$Z$ is smooth of dimension $n$, then Poincar{\'e} duality implies that
 $$
 \chi^c_{p,q}(Z)= \chi_{n-p,n-q}(Z).
 $$

\begin{definition}
For \textit{any} complex algebraic variety $Z$ (not necessarily
smooth, compact or irreducible), we define  its {\em Hodge
polynomial} \cite{DK} as
 $$
 e(Z)=e(Z)(u,v)= \sum_{p,q} (-1)^{p+q}\chi^c_{p,q}(Z) u^p v^q\, .
 $$
\end{definition}

Note that if $Z$ is smooth and projective then the mixed Hodge
structure on $H_c^k(Z)$ is pure of weight $k$. This means that
$\Gr_k^W H_c^k(Z) =H_c^k(Z)=H^k(Z)$ and the other pieces $\Gr_m^W
H_c^k(Z)=0$, $m\neq k$. So
   $$
   \chi_{p,q}(Z) = \chi^c_{p,q}(Z)=(-1)^{p+q} h^{p,q}(Z),
   $$
where $h^{p,q}(Z)$ is the usual Hodge number of $Z$. In this case,
  $$
  e(Z)(u,v)= \sum_{p,q} h^{p,q}(Z) u^p v^q\,
  $$
is the (usual) Hodge polynomial of $Z$. Note that in this case,
the Poincar{\'e} polynomial of $Z$ is
  $$
  P_t(Z)=\sum_k b^k(Z) t^k= \sum_k \left(
  \sum_{p+q=k} h^{p,q}(Z) \right) t^k= e(Z)(t,t).
  $$
where $b^k(Z)$ is the $k$-th Betti number of $Z$.

The following result is a slight extension of that of \cite{Du}.

\begin{theorem} \label{thm:Du}
 Let $Z$ be a complex algebraic variety. Suppose that $Z$ is
 a finite disjoint union $Z=Z_1\sqcup \cdots \sqcup Z_n$, where the
 $Z_i$ are algebraic subvarieties. Then
  $$
  \chi^c_{p,q}(Z)= \sum_i \chi^c_{p,q} (Z_i).
  $$
 Equivalently, $e(Z)=\sum_i e(Z_i)$.
\end{theorem}

\begin{dem}
It is enough to see that for a complex quasi-projective variety
$Z$, a closed sub\-va\-rie\-ty $Y$ of $Z$ and $U=Z-Y$, we have
  $$
  e(Z) = e(Y) + e(U).
  $$

To see this, consider the long exact sequence
 $$
 \ldots \to H_c^k(U) \to H_c^k(Z)\to H_c^k(Y) \to H_c^{k+1}(U)\to
 \ldots
 $$
of cohomology groups with compact supports. The maps in this exact
sequence are com\-pa\-ti\-ble with the weight and Hodge
filtrations. The induced sequence on the $(p,q)$-pieces of the
Hodge decomposition remains exact for all $p,q$:
 $$
 \ldots \to \Gr_F^p \Gr_{p+q}^W H_c^k(U) \to \Gr_F^p \Gr_{p+q}^W
 H_c^k(Z)\to
 \Gr_F^p \Gr_{p+q}^W H_c^k(Y) \to \Gr_F^p \Gr_{p+q}^W
 H_c^{k+1}(U)\to
 \ldots
 $$
This means that the Euler characteristics $\chi^c_{p,q}$ satisfy
  $$
   \chi_{p,q}^c(Z) = \chi^c_{p,q}(Y) + \chi^c_{p,q}(U),
  $$
from where $e(Z) = e(Y) + e(U)$ follows.
\end{dem}

Let us give a collection of simple examples:
\begin{itemize}
 \item Let $Z=\PP^{n}$, then $e(Z)= 1+uv+(uv)^2+\cdots +
 (uv)^{n}=(1-(uv)^{n+1})/(1-uv)$. For future reference,
 we shall denote
   \begin{equation}\label{eqn:Pn}
    e_n := e
    (\PP^{n-1})=e
    (\PP(\CC^n))
   =\frac{1-(uv)^{n}}{1-uv}\ .
   \end{equation}
 \item Let $\Jac^d X$ be the Jacobian of (any) degree $d$ of a
 (smooth, projective) complex curve $X$ of genus $g$. Then
  \begin{equation}\label{eqn:Jac}
  e(\Jac^d X)=(1+u)^g(1+v)^g.
   \end{equation}
 \item Let $X$ be a curve (smooth, projective) complex curve of genus
 $g$, and $k\geq 1$. The Hodge polynomial of the symmetric product
  $\Sym^k X$ is computed in \cite{Bur},
  \begin{equation}\label{eqn:sym}
  e(\Sym^k X)=
 \coeff_{x^0}\frac{(1+ux)^{g}(1+vx)^{g}}{(1-x)(1-uvx)x^{k}}\, .
  \end{equation}
  \end{itemize}

\medskip

\begin{lemma}\label{lem:vb}
  Suppose that  $\pi:Z\to Y$ is an algebraic fiber bundle with fiber $F$ which is
 locally trivial in the Zariski topology, then $e(Z)=e(F)\,e(Y)$.
 (In particular this is true for $Z=F\times Y$.)
\end{lemma}

\begin{dem}
First we need to see that $e(F\times Y)=e(F)\, e(Y)$. Using
Theorem \ref{thm:Du}, we see that it is enough to see it for
quasi-projective smooth varieties $F$ and $Y$. We may prove it by
induction on the dimension of $F\times Y$. First, if both $F$ and
$Y$ are smooth and projective, then $e(F\times Y)=e(F)\, e(Y)$ by
usual Hodge theory. Second, assume that $F$ is smooth projective
and $Y$ is smooth quasi-projective. Then write $Y=Y_1-Y_2$, where
$Y_1$ is a smooth projective variety and $Y_2\subset Y_1$ is a
smaller dimensional subvariety (this is possible by Hironaka's
result on embedded resolution of singularities). Then $e(F\times
Y)=e(F\times Y_1) - e(F\times Y_2)= e(F) e(Y_1)-e(F) e(Y_2)=e(F)
(e(Y_1)-e(Y_2))=e(F)e(Y)$. The general case, when $F$ and $Y$ are
both smooth quasi-projective, now follows from the previous case
by writing $F=F_1-F_2$, with $F_1$ smooth projective variety and
$F_2\subset F_1$ a smaller dimensional subvariety.

Decompose $Y=Y_1\sqcup \cdots \sqcup Y_n$ into a finite union of
disjoint locally closed subvarieties such that $Z_i=\pi^{-1}(Y_i)
\cong F \times Y_i$ (using the  local triviality in the Zariski
topology). Then $e(Z)=\sum e(Z_i)=\sum e(F)e(Y_i)=e(F)e(Y)$.
\end{dem}

\begin{remark}
Consider the abelian group $\cV{}ar$ generated by all algebraic
varieties subject to the following equivalence relation: if
$Z=Z_1\sqcup Z_2$, then set $[Z] = [Z_1] +[Z_2]$. The cartesian
product of varieties yields an algebra structure for $\cV{}ar$.
For instance, for a (locally trivial in the Zariski topology)
bundle $Z\to Y$ with fiber $F$, we have $[Z]=[F]\cdot [Y]$. The
Hodge polynomial $[Z]\mapsto e(Z)$ gives a morphism (of algebras)
from $\cV{}ar$ to $\ZZ[u,v]$.
\end{remark}

%%%%%%%%%%%%%%%%%%%%%%%%%%%%%%%%%%%%%%%%%%%%%%%%%%%%%%%%%%%%%%%%%%%%%%%
\section{Moduli spaces of triples}
\label{sec:stable-triples}
%%%%%%%%%%%%%%%%%%%%%%%%%%%%%%%%%%%%%%%%%%%%%%%%%%%%%%%%%%%%%%%%%%%%%%%

Let $X$ be a smooth projective curve over the complex numbers of
genus $g\geq 2$. A \emph{holomorphic triple} $T =
(E_{1},E_{2},\phi)$ on $X$ consists of two holomorphic vector
bundles $E_{1}$ and $E_{2}$ over $X$, of ranks $n_1$ and $n_2$ and
degrees $d_1$ and $d_2$, respectively, and a holomorphic map $\phi
\colon E_{2} \to E_{1}$. We refer to $(n_1,n_2,d_1,d_2)$ as the
\emph{type} of $T$, to $(n_1,n_2)$ as the \emph{rank} of $T$, and
to $(d_1,d_2)$ as the \emph{degree} of $T$.

A homomorphism from $T' = (E_1',E_2',\phi')$ to $T = (E_1,E_2,\phi)$ is a
commutative diagram
  \begin{displaymath}
  \begin{CD}
    E_2' @>\phi'>> E_1' \\
    @VVV @VVV  \\
    E_2 @>\phi>> E_1,
  \end{CD}
  \end{displaymath}
where the vertical arrows are homomorphisms. A triple
$T'=(E_1',E_2',\phi')$ is a subtriple of $T = (E_1,E_2,\phi)$ if
$E_1'\subset E_1$ and $E_2'\subset E_2$ are subbundles,
$\phi(E_2')\subset E_1'$ and $\phi'=\phi|_{E_2'}$. A subtriple
$T'\subset T$ is called \emph{proper} if $T'\neq 0 $ and $T'\neq
T$. The quotient triple $T''=T/T'$ is given by $E_1''=E_1/E_1'$,
$E_2''=E_2/E_2'$ and $\phi'' \colon E_2''\to E_1''$ being the map
induced by $\phi$. We usually denote by $(n_1',n_2',d_1',d_2')$
and $(n_1'',n_2'',d_1'',d_2'')$, the types of the subtriple $T'$
and the quotient triple $T''$.

\begin{definition} \label{def:s-slope}
For any $\s \in \RR$ the \emph{$\s$-slope} of $T$ is defined by
 $$
   \mu_{\s}(T)  =
   \frac{d_1+d_2}{n_1+n_2} + \s \frac{n_{2}}{n_{1}+n_{2}}\ .
 $$
To shorten the notation, we define the \emph{$\mu$-slope} and
\emph{$\lambda$-slope} of the triple $T$ as $\mu=\mu(E_{1} \oplus E_{2})=
\frac{d_1+d_2}{n_1+n_2}$ and $\lambda=\frac{n_{2}}{n_{1}+n_{2}}$, so that
$\mu_{\s}(T)=\mu+\s \lambda$.
\end{definition}

\begin{remark} \label{rem:slopes}
For any triple $T=(E_1,E_2,\phi)$, if $T'=(E_1',E_2',\phi')\subset T$ is a
proper subtriple and $T''=T/T'$ is the corresponding quotient, then only
one of the following situations is possible:
  \begin{enumerate}
  \item[(1)] $\mu_\s(T)=\mu_\s(T')=\mu_\s(T'')$\ ;
  \item[(2)] $\mu_\s(T')<\mu_\s(T)<\mu_\s(T'')$\ ;
  \item[(3)] $\mu_\s(T')>\mu_\s(T)>\mu_\s(T'')$\ .
  \end{enumerate}
In fact, letting $t=\frac{n_1'+n_2'}{n_1+n_2}\in (0,1)$, the next equality
holds
  \begin{equation*}
  \mu_\s(T)=t\cdot\mu_\s(T')+(1-t)\cdot\mu_\s(T'')
  \end{equation*}
Also note the following equality of $\lambda$-slopes
  \begin{equation}\label{eqn:lambda}
  \lambda=t\cdot\lambda' +(1-t)\cdot\lambda''\, ,
  \end{equation}
where $\lambda$, $\lambda'$ and $\lambda''$ are the $\lambda$-slopes of
$T$, $T'$ and $T''$, respectively.
\end{remark}

\begin{definition}\label{def:sigma-stable}
We say that a triple $T = (E_{1},E_{2},\phi)$ is \emph{$\s$-stable} if
  $$
  \mu_{\s}(T') < \mu_{\s}(T) \quad
  \mbox{ (equivalently $\mu_\s(T/T')>\mu_\s(T)\,$)},
  $$
for any proper subtriple $T' = (E_{1}',E_{2}',\phi')$. We define
\emph{$\s$-semistability} by replacing the above strict inequality with a
weak inequality. A triple is called \emph{$\s$-polystable} if it is the
direct sum of $\s$-stable triples of the same $\s$-slope. It is
\emph{$\s$-unstable} if it is not $\s$-semistable, and \emph{strictly
$\s$-semistable} if it is $\s$-semistable but not $\s$-stable. A
$\s$-destabilizing subtriple $T'\subset T$ is a proper subtriple
satisfying $\mu_{\s}(T') \geq \mu_{\s}(T)$.
\end{definition}

We denote by
  $$
  \cN_\s = \cN_\s(n_1,n_2,d_1,d_2)
  $$
the moduli space of $\s$-polystable triples $T =
(E_{1},E_{2},\phi)$ of type $(n_1,n_2,d_1,d_2)$, and drop the type
from the notation when it is clear from the context.  This moduli
space is constructed in \cite{BGP} by using dimensional reduction.
A direct construction is given by Schmitt \cite{Sch} using
geometric invariant theory. The subspace of $\s$-stable triples is
denoted by $ \cN_\s^s = \cN_\s^s(n_1,n_2,d_1,d_2)$.

Given a triple  $T=(E_1,E_2,\phi)$ one has the dual triple
$T^*=(E_2^*,E_1^*,\phi^*)$, where $E_i^*$ is the dual of $E_i$ and
$\phi^*$ is the transpose of $\phi$. The following  is not difficult to
prove.

\begin{proposition}\cite[Proposition 3.16]{BGP}
\label{prop:duality} The  $\s$-(semi)stability of $T$   is equivalent to
the $\s$-(semi)stability of $T^*$. The map $T\mapsto T^*$ defines an
isomorphism
 $$
 \cN_\s(n_1,n_2,d_1,d_2) \cong \cN_\s(n_2,n_1,-d_2,-d_1)\,
 .
 $$ \hfill $\Box$
\end{proposition}

This can be used to restrict our study to $n_1\geq n_2$ and appeal to
duality for $n_1<n_2$.

\medskip

The particular case where $n_1=0$ or $n_2=0$ reduces to the case of
bundles. Let $M(n,d)$ denote the moduli space of polystable vector bundles
of rank $n$ and degree $d$ over $X$. This moduli space is projective. We
also denote by $M^s(n,d)$ the open subset of stable bundles, which is
smooth of dimension $n^2(g-1)+1$. If $\GCD(n,d)=1$, then
$M(n,d)=M^s(n,d)$.

\begin{lemma} \label{lem:(n,0)}
There are isomorphisms $\cN_\s(n,0,d,0)\cong M(n,d)$ and
$\cN_\s^s(n,0,d,0)\cong M^s(n,d)$, for all $\s\in \RR$. In particular,
$\cN_\s(1,0,d,0)= \cN_\s^s(1,0,d,0) \iso \Jac^{d}X$, for any $\s\in \RR$.
%\hfill $\Box$
\end{lemma}

There are certain necessary conditions in order for $\s$-semistable
triples to exist. Let $\mu_i=\mu(E_i)=d_i/n_i$ stand for the slope of
$E_i$, for $i=1,2$. We write
  \begin{align*}
  \s_m = &\, \mu_1-\mu_2\ ,  \\
  \s_M = & \left(1+ \frac{n_1+n_2}{|n_1 - n_2|}\right)(\mu_1 -
\mu_2)\ ,
      \qquad \mbox{if $n_1\neq n_2$\ .}
  \end{align*}

\begin{proposition}\cite{BGPG} \label{prop:alpha-range}
The moduli space $\cN_\s(n_1,n_2,d_1,d_2)$ is a complex projective
variety. A necessary condition for $\cN_\s(n_1,n_2,d_1,d_2)$ to be
non-empty is
  $$
  \begin{array}{ll}
   0\leq \s_m \leq \s \leq \s_M,
  \quad &\text{if }\ n_1\neq n_2\ , \ n_1\neq 0\ ,\  n_2\neq 0
\ ,\\
   0\leq \s_m \leq \s, \quad &\text{if }\ n_1= n_2\neq 0\ .
  \end{array}
  $$%\hfill $\Box$
\end{proposition}

We shall denote by $I\subset\RR$ the following interval
  $$
  I=\left\{ \begin{array}{ll} [\s_m,\s_M], \qquad &\hbox{if
  $n_1\neq n_2$, $n_1\neq 0$, $n_2\neq 0$,} \\
  {} [\s_m,\infty), & \hbox{if $n_1=n_2\neq 0$,} \\
  \RR, & \hbox{if $n_1=0$ or $n_2=0$.} \end{array} \right.
  $$

\bigskip

To study the dependence of the moduli spaces $\cN_\s$ on the parameter, we
need to introduce the concept of critical value.

\begin{proposition}\label{prop:sigma-coeff}
 Let $\s_0 \in I$ and let $T=(E_1,E_2,\phi) \in \cN_{\s_0}
 (n_1,n_2,d_1,d_2)$ be a strictly $\s_0$-semistable triple. Then
 one of the following conditions holds:
 \begin{enumerate}
  \item[(1)] For all $\s_0$-destabilizing subtriples
$T'=(E_1',E_2',\phi')$,
  we have $\displaystyle\lambda'=\frac{n_2'}{n_1'+n_2'}=
  \lambda=\frac{n_2}{n_1+n_2}$. Then
  $T$ is strictly $\s$-semistable for $\s\in
  (\s_0-\varepsilon, \s_0+\varepsilon)$, for some small
  $\varepsilon>0$.

  \item[(2)] There exists a $\s_0$-destabilizing subtriple
$T'=(E_1',E_2',\phi')$
  with $\displaystyle\lambda'=\frac{n_2'}{n_1'+n_2'} \neq
 \lambda=\frac{n_2}{n_1+n_2}$. Then either: \smallskip
  \begin{itemize}
  \item[(i)] $\lambda'>\lambda$. Then for
  any $\s>\s_0$, $T$~is $\s$-unstable.
  \item[(ii)] $\lambda'<\lambda$.
  Then for any $\s<\s_0$, $T$~is $\s$-unstable.
  \end{itemize}
 \end{enumerate}

 Conversely, if $T$ is a triple such that it is
 $\s$-semistable for $\s<\s_0$ and $\s$-unstable
 for $\s>\s_0$ (resp.\ $\s$-semistable
 for $\s>\s_0$ and $\s$-unstable
 for $\s<\s_0$) then $T$ is strictly $\s_0$-semistable
 and Case {\rm (i)} (resp.\ Case {\rm (ii)}) holds.
\end{proposition}

\begin{dem}
If $T'$ is a $\s_0$-destabilizing subtriple, then $\mu_{\s_0}(T')=
\mu_{\s_0}(T)$, i.e., $\mu'+\s_0 \lambda'=\mu +\s_0 \lambda$. If (2)
happens, then $\lambda'\neq \lambda$ (which is equivalent to $n_1'n_2\neq
n_1n_2'$). We have the following pictures of the $\s$-slopes of $T$ and
$T'$, as functions of $\s$.

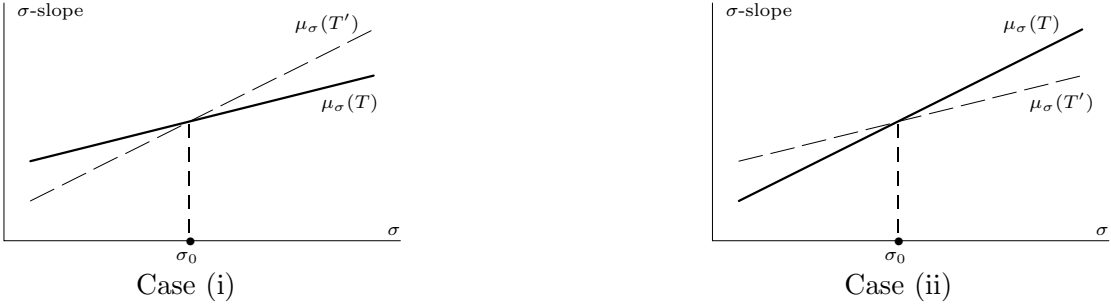
\begin{minipage}{160pt}
\begin{picture}(150,140)(0,-20)
    \put(150,0){\path(0,0)(-150,0)(-150,90)}
    \multiput(10,15)(12,6){11}{\line(2,1){9}}%0}}
    \thicklines
    \put(10,30){\line(4,1){130}}
    \thinlines
    \put(110,80){\tiny$\mu_\s(T')$}
    \put(120,50){\tiny$\mu_\s(T)$}
    \multiput(70,3)(0,9){5}{\line(0,1){5}}
    \put(68.5,-2){\tiny$\bullet$}
    \put(65,-7.5){\tiny$\s_0$}
    \put(50,-20){Case (i)}
    \put(145,2){\tiny$\s$}
    \put(5,85){\tiny$\s$-slope}
\end{picture}
\end{minipage}\hfill
\begin{minipage}{160pt}
\begin{picture}(150,140)(0,-20)
    \put(150,0){\path(0,0)(-150,0)(-150,90)}
    \multiput(10,30)(12,3){11}{\line(4,1){9}}%0}}
    \thicklines
    \put(10,15){\line(2,1){130}}
    \thinlines
    \put(110,80){\tiny$\mu_\s(T)$}
    \put(120,50){\tiny$\mu_\s(T')$}
    \multiput(70,3)(0,9){5}{\line(0,1){5}}
    \put(68.5,-2){\tiny$\bullet$}
    \put(65,-7.5){\tiny$\s_0$}
    \put(50,-20){Case (ii)}
    \put(145,2){\tiny$\s$}
    \put(5,85){\tiny$\s$-slope}
\end{picture}
\end{minipage}

\medskip
Wherever the graph of $\mu_\s(T')$ is above that of $\mu_\s(T)$,
the triple $T$ is $\s$-unstable. This yields Cases (i) and (ii).

On the other hand, if $T'$ is a $\s_0$-destabilizing subtriple and
Case (1) happens, then the equality $\mu'+\s_0 \lambda'=\mu +\s_0
\lambda$ and $\lambda'=\lambda$ imply $\mu'=\mu$, hence
$\mu_{\s}(T')= \mu_{\s}(T)$, for all $\s$, so $T'$ is
$\s$-destabilizing. But for any other subtriple $\tilde T'\subset
T$, either $\mu_{\s_0}(\tilde T')<\mu_{\s_0}(T)$ in which case
$\mu_{\s}(\tilde T')<\mu_{\s}(T)$ for $\s$ close to $\s_0$, or
$\mu_{\s_0}(\tilde T')=\mu_{\s_0}(T)$ in which case
$\mu_{\s}(\tilde T')=\mu_{\s}(T)$ by the above argument. So $T$ is
strictly $\s$-semistable for $\s$ close to $\s_0$.

The converse statement follows easily from the pictures of the $\s$-slopes
of $T$ and of a subtriple $T'$ with $\mu_\s(T')>\mu_\s(T)$ for $\s>\s_0$
(resp.\ for $\s<\s_0$).
\end{dem}

\begin{remark}\label{rem:sigma-independent}
 We call the phenomenon in
 Case (1) in Proposition
 \ref{prop:sigma-coeff} \ \emph{$\s$-independent semistability}.
 Note that this implies simultaneously that
  $$
  \frac{n'_2}{n_1'+n_2'} = \frac{n_2}{n_1+n_2}\ ,\;\; \mbox{and}
\;
  \;\;\; \frac{d'_1+d'_2}{n_1'+n_2'} = \frac{d_1+d_2}{n_1+n_2}\
.
  $$
 Hence this cannot happen if $\GCD(n_1,n_2,d_1+d_2) = 1$.
\end{remark}

\begin{definition}\label{def:critical}
The values of $\s\in I$ for which Case (2) in Proposition
\ref{prop:sigma-coeff} occurs are called \emph{critical values}. If
$\s_c\in I$ is a critical value then there exist integers $n'_1$, $n'_2$,
$d'_1$ and $d'_2$ such that
 $$
  \frac{d'_1+d'_2}{n'_1+n'_2} + \s_c\frac{n'_2}{n'_1+n'_2}
  = \frac{d_1+d_2}{n_1+n_2} + \s_c\frac{n_2}{n_1+n_2},
 $$
equivalently,
 \begin{equation}\label{eqn:sigmac}
 \s_c=\frac{(n_1+n_2)(d_1'+d_2')-(n_1'+n_2')(d_1+d_2)}{n_1'n_2-n_1n_2'},
 \end{equation}
with $0 \le n'_i \leq n_i$ and  $n_1'n_2\neq n_1n_2'$. We call the
numbers (\ref{eqn:sigmac}) \textit{virtual critical values}. In
general, not all the virtual critical values are critical values.
We say that $\s\in I$ is \emph{generic} if it is not critical.
\end{definition}

\begin{remark}\label{rem:directsum}
In the Case (2) of Proposition \ref{prop:sigma-coeff}, the
quotient triple $T''=T/T'$ satisfies
$\mu_{\s_c}(T'')=\mu_{\s_c}(T)$ by Remark \ref{rem:slopes}. By
(\ref{eqn:lambda}) the graph of $\mu_\s(T'')$ is a line which goes
on the opposite side of $\mu_\s(T')$, i.e., it has smaller slope
than $\lambda$ in Case (i), and bigger slope than $\lambda$ in
Case (ii).

In particular, for a $\s_c$-semistable decomposable triple
$T=T'\oplus T''$, with $\mu_{\s_c}(T)= \mu_{\s_c}(T')=
\mu_{\s_c}(T'')$ and $n_1'n_2\neq n_1n_2'$, $T$ is $\s$-unstable
for any $\s\neq\s_c$. For instance, if $T=(E_1,E_2,\phi)$ has
$\phi=0$ (and $n_1n_2\neq 0$), then $T=(E_1,0,0)\oplus (0,E_2,0)$,
so it is not $\s$-stable for any value of $\s$ (cf.\ \cite[Lemma
3.5]{BGP}).
\end{remark}

\begin{lemma}\label{lem:both-i-and-ii}
Given a triple $T$, let $I_s(T)=\{\s\in \RR \, |\, \hbox{$T$ is
$\s$-stable}\}$ and $I_{ss}(T)=\{\s\in \RR \, |\, \hbox{$T$ is
$\s$-semistable}\}$. Then
 \begin{itemize}
  \item[(a)] $I_{ss}(T)$ is one of the
  following: $\emptyset$, $\{\s_c\}$, $[\s_c,\s_c']$,
  $[\s_c, \infty)$, $(-\infty,\s_c]$, $\RR$. (Here $\s_c$, $\s_c'$
  are virtual critical values.)
  \item[(b)] $I_s(T)$ is either empty or it is the interior of
  $I_{ss}(T)$.
 \end{itemize}
\end{lemma}

\begin{dem}
For each subtriple $T'\subset T$, let $I_{ss}(T)_{T'}=\{\s\in \RR \, |\,
\mu_\s(T')\leq \mu_\s(T)\}$ and $I_{s}(T)_{T'}=\{\s\in \RR \,
|\,\mu_\s(T')<\mu_\s(T)\}$. Then we have the following possibilities:
 \begin{itemize}
 \item If $\lambda'>\lambda$ then $I_{ss}(T)_{T'}=(-\infty,\s_c]$
 and $I_{s}(T)_{T'}=(-\infty, \s_c)$.
 \item If $\lambda'<\lambda$ then $I_{ss}(T)_{T'}=[\s_c,\infty)$
 and $I_{s}(T)_{T'}=(\s_c,\infty)$.
 \item If $\lambda'=\lambda$ and $\mu'<\mu$ then
$I_{ss}(T)_{T'}=\RR$ and
 $I_{s}(T)_{T'}=\RR$.
  \item If $\lambda'=\lambda$ and $\mu'=\mu$ then
$I_{ss}(T)_{T'}=\RR$ and
 $I_{s}(T)_{T'}=\emptyset$.
  \item If $\lambda'=\lambda$ and $\mu'>\mu$ then
$I_{ss}(T)_{T'}=\emptyset$ and
 $I_{s}(T)_{T'}=\emptyset$.
 \end{itemize}
where $\s_c$ is always a virtual critical value. The set of possible
values (\ref{eqn:sigmac}) for $\s_c$ is a discrete set of the real line.
So the result follows easily from $I_s(T)=\bigcap I_s(T)_{T'}$ and
$I_{ss}(T)=\bigcap I_{ss}(T)_{T'}$, where the intersection runs over the
subtriples $T'$ of $T$.
\end{dem}

\begin{remark}\label{rem:otromas}
 If $n_1=0$ or $n_2=0$ then $I_s(T)$ and $I_{ss}(T)$ are either
 $\emptyset$ or $\RR$. If $n_1\neq 0$, $n_2\neq 0$ and $n_1\neq
 n_2$ then the last three cases of (a) of Lemma
 \ref{lem:both-i-and-ii} do not happen since $I_{ss}(T)\subset
 I=[\s_m,\s_M]$. If $n_1=n_2\neq 0$ then the last two cases of (a) of
Lemma
 \ref{lem:both-i-and-ii} do not happen since $I_{ss}(T)\subset
 I=[\s_m,\infty)$.
\end{remark}

\begin{proposition} \label{prop:triples-critical-range}
Fix $(n_1,n_2,d_1,d_2)$. Then
 \begin{enumerate}
 \item[(1)] The critical values are a finite number of values $\s_c
\in I$.
 \item[(2)] The stability and semistability criteria  for two values of
$\s$
  lying between two con\-se\-cu\-ti\-ve critical values are equivalent; thus
  the corresponding moduli spaces are isomorphic.
 \item[(3)] If $\s$ is generic and $\GCD(n_1,n_2,d_1+d_2) = 1$,
  then $\s$-semistability is equivalent to $\s$-stability, i.e.,\
  $\cN_\s=\cN_\s^s$.
 \item[(4)] If $\GCD(n_1, n_2,d_1+d_2)=1$ then the moduli spaces
$\cN_\s^s$
 are fine moduli spaces, i.e.,\ there is a universal triple $\cT
 \to X\times\cN_{\s}^s$.
 \end{enumerate}
\end{proposition}

\begin{dem}
 Items (1), (2) and (3) follow from \cite[Proposition 2.6]{BGPG}.
 Item (4) is in \cite{Sch}.
\end{dem}

\begin{remark}\label{rem:duality}
 Under the isomorphism of Proposition \ref{prop:duality}, the
 critical values for triples of type $(n_1,n_2,d_1,d_2)$ and
 of type $(n_2,n_1,-d_2,-d_1)$ correspond. This can be seen at the
 level of virtual critical values, since the formula
 (\ref{eqn:sigmac}) is unchanged by the transformation $(n_1,n_2,d_1,d_2,
 n_1',n_2',d_1',d_2')\mapsto
 (n_2,n_1,-d_2,-d_1,n_2-n_2',n_1-n_1',d_2'-d_2,d_1'-d_1)$.
\end{remark}

\begin{remark}\label{rem:critical(n,0)}
There are no critical values for triples of rank $(n,0)$ or $(0,n)$.
\end{remark}

The moduli space of triples of rank $(1,1)$ is easily described as
follows:

\begin{lemma}\label{lem:type(1,1)}
 \begin{itemize}
  \item The moduli space
  $\cN_\s(1,1,d_1,d_2)$ is empty if $d_1<d_2$.
  \item If $d_1\geq
  d_2$ then $\cN_\s(1,1,d_1,d_2) = \cN^s_\s(1,1,d_1,d_2)
  \iso \Jac^{d_2}X \x \Sym^{d_1-d_2}X$, for any
  $\s>\s_m=d_1-d_2$; it is empty if $\s<\s_m$; and
  $\cN_{\s_m}=\Jac^{d_1}X\times \Jac^{d_2}X$,
  $\cN_{\s_m}^s=\emptyset$.
 \end{itemize}
\end{lemma}

\begin{dem}
  The first item follows from Proposition \ref{prop:alpha-range}
(otherwise,
  there are no non-trivial maps $E_2\to E_1$ and by Remark
\ref{rem:directsum}
  such triples cannot be $\sigma$-stable for any $\sigma$).
  For the second item, the moduli space is parametrizing
  triples of the form $\phi:E_2\to E_1$, where both $E_1$ and $E_2$
  are line bundles. The only possible proper subtriple is $0\to E_1$
  and this violates $\s$-stability for $\s\leq \s_m$. If
  $\s=\s_m$, then the triple is strictly $\s_m$-semistable and, being
  polystable, it is $(0,E_2,0)\oplus (E_1,0,0)$. (Another
consequence
  is that $\s_m$ is the only critical value.) If
  $\s>\s_m$ then it cannot happen that $\phi=0$. So the triples
$T=(E_1,E_2,\phi)$ are
  parametrized by $(E_2,Z(\phi))\in \Jac^{d_2}X \x \Sym^{d_1-d_2}X$.
\end{dem}
%\mbox{}\vskip-3\baselineskip\mbox{}

%%%%%%%%%%%%%%%%%%%%%%%%%%%%%%%%%%%%%%%%%%%%%%%%%%%
\section{Extensions of triples}
\label{sec:extensions-of-triples}
%%%%%%%%%%%%%%%%%%%%%%%%%%%%%%%%%%%%%%%%%%%%%%%%%%%

The homological algebra of triples is controlled by the hypercohomology of
a certain complex of sheaves which appears when studying infinitesimal
deformations \cite[Section 3]{BGPG}. Let $T'=(E'_1,E'_2,\phi')$ and
$T''=(E''_1,E''_2,\phi'')$ be two triples of types
$(n_{1}',n_{2}',d_{1}',d_{2}')$ and $(n_{1}'',n_{2}'',d_{1}'',d_{2}'')$,
respectively. Let $\Hom(T'',T')$ denote the linear space of homomorphisms
from $T''$ to $T'$, and let $\Ext^1(T'',T')$  denote the linear space of
equivalence classes of extensions of the form
 $$
  0 \lto T' \lto T \lto T'' \lto 0,
 $$
where by this we mean a commutative  diagram
  $$
  \begin{CD}
  0@>>>E_1'@>>>E_1@>>> E_1''@>>>0\\
  @.@A\phi' AA@A \phi AA@A \phi'' AA\\
  0@>>>E'_2@>>>E_2@>>>E_2''@>>>0.
  \end{CD}
  $$
To analyze $\Ext^1(T'',T')$ one considers the complex of sheaves
 \begin{equation} \label{eqn:extension-complex}
    C^{\bullet}(T'',T') \colon ({E_{1}''}^{*} \otimes E_{1}')
\oplus
  ({E_{2}''}^{*} \otimes E_{2}')
  \overset{c}{\lto}
  {E_{2}''}^{*} \otimes E_{1}',
 \end{equation}
where the map $c$ is defined by
 $$
 c(\psi_{1},\psi_{2}) = \phi'\psi_{2} - \psi_{1}\phi''.
 $$

\begin{proposition}[{\cite[Proposition 3.1]{BGPG}}]
  \label{prop:hyper-equals-hom}
  There are natural isomorphisms
  \begin{align*}
    \Hom(T'',T') &\cong \HH^{0}(C^{\bullet}(T'',T')), \\
    \Ext^{1}(T'',T') &\cong \HH^{1}(C^{\bullet}(T'',T')),
  \end{align*}
and a long exact sequence associated to the complex $C^{\bullet}(T'',T')$:
 $$
 \begin{array}{c@{\,}c@{\,}c@{\,}l@{\,}c@{\,}c@{\,}c}
  0 &\lto \mathbb{H}^0(C^{\bullet}(T'',T')) &
  \lto & H^0(({E_{1}''}^{*} \otimes E_{1}') \oplus ({E_{2}''}^{*}
\otimes
  E_{2}'))
  & \lto &  H^0({E_{2}''}^{*} \otimes E_{1}') \\
    &  \lto \mathbb{H}^1(C^{\bullet}(T'',T')) &
  \lto &  H^1(({E_{1}''}^{*} \otimes E_{1}') \oplus
({E_{2}''}^{*} \otimes
  E_{2}'))
 &  \lto & H^1({E_{2}''}^{*} \otimes E_{1}') \\
 &   \lto \mathbb{H}^2(C^{\bullet}(T'',T')) & \lto & 0. & &
 \end{array}
 $$ %\hfill $\Box$
\end{proposition}

We introduce the following notation:
\begin{align*}
  h^{i}(T'',T') &= \dim\HH^{i}(C^{\bullet}(T'',T')), \\
%\notag \\
  \chi(T'',T') &= h^0(T'',T') - h^1(T'',T') + h^2(T'',T').
\end{align*}

\begin{proposition}[{\cite[Proposition 3.2]{BGPG}}]
  \label{prop:chi(T'',T')}
  For any holomorphic triples $T'$ and $T''$ we have
  \begin{align*}
    \chi(T'',T') &= \chi({E_{1}''}^{*} \otimes E_{1}')
    + \chi({E_{2}''}^{*} \otimes E_{2}')
    - \chi({E_{2}''}^{*} \otimes E_{1}')  \\
    &= (1-g)(n''_1 n'_1 + n''_2 n'_2 - n''_2 n'_1) + n''_1
d'_1 - n'_1 d''_1
    + n''_2 d'_2 - n'_2 d''_2
    - n''_2 d'_1 + n'_1 d''_2,
  \end{align*}
where $\chi(E)=\dim H^0(E) - \dim H^1(E)$ is the Euler characteristic of
$E$. %\hfill $\Box$
\end{proposition}

\begin{proposition}[{\cite[Proposition 3.5]{BGPG}}]
\label{prop:h0-vanishing}
   Suppose that $T'$ and $T''$ are $\s$-semistable, for some
   value of $\s$.
 \begin{enumerate}
  \item[(1)] If $\mu_\s(T')<\mu_\s (T'')$ then
  $\HH^{0}(C^{\bullet}(T'',T')) = 0$.
  \item[(2)] If $\mu_\s(T')=\mu_\s (T'')$ and $T''$ is
   $\s$-stable, then
  $$
     \HH^{0}(C^{\bullet}(T'',T')) \cong
     \begin{cases}
       \CC \quad &\text{if $T' \cong T''$} \\
       0 \quad &\text{if $T' \not\cong T''$}.
     \end{cases}
  $$
 \end{enumerate} %\hfill $\Box$
\end{proposition}

Since the  space of infinitesimal deformations of $T$ is isomorphic to
$\HH^{1}(C^{\bullet}(T,T))$, the previous results also apply to studying
deformations of a holomorphic triple $T$.

\begin{theorem}[{\cite[Theorem 3.8]{BGPG}}]\label{thm:smoothdim}
Let $T=(E_1,E_2,\phi)$ be an $\s$-stable triple of type
$(n_1,n_2,d_1,d_2)$.
 \begin{enumerate}
 \item[(1)] The Zariski tangent space at the point defined by $T$
 in the moduli space of stable triples  is isomorphic to
 $\HH^{1}(C^{\bullet}(T,T))$.
 \item[(2)] If\/ $\HH^{2}(C^{\bullet}(T,T))= 0$, then the moduli
space of
 $\s$-stable triples is smooth in  a neighbourhood of the point
 defined by $T$.
% \item[(3)] $\HH^{2}(C^{\bullet}(T,T))= 0$ if and only if the
%homomorphism
% $$
%  H^1((E_1^* \otimes E_1) \oplus (E_2^* \otimes E_2))
%  \lto  H^1(E_2^* \otimes E_1)
% $$
% in the corresponding long exact sequence is surjective.
 \item[(3)] At a smooth point $T\in \cN^s_\s(n_1,n_2,d_1,d_2)$ the
 dimension of the moduli space of $\s$-stable triples is
 \begin{align*}
  \dim \cN^s_\s(n_1,n_2,d_1,d_2)
  &= h^{1}(T,T) = 1 - \chi(T,T) \\ %\notag \\
  &= (g-1)(n_1^2 + n_2^2 - n_1 n_2) - n_1 d_2 + n_2 d_1 + 1.
 \end{align*}
 \item[(4)] Let $T=(E_1,E_2,\phi)$ be a $\s$-stable triple. If $T$
 is injective or surjective (meaning that $\phi:E_2\to E_1$ is
 injective or surjective) then the moduli space is smooth at $T$.
 \end{enumerate}
\end{theorem}

\subsection*{Description of the flip loci}
Fix the type $(n_1,n_2,d_1,d_2)$ for the moduli spaces of
holomorphic triples. Now we want to describe the differences
between two spaces $\cN^s_{\s_1}$ and $\cN^s_{\s_2}$ when $\s_1$
and $\s_2$ are separated by a critical value. Let $\s_c\in I$ be a
critical value and set
 $$
 \scp = \s_c + \varepsilon,\quad \scm = \s_c -
 \varepsilon,
 $$
where $\varepsilon > 0$ is small enough so that $\s_c$ is the only critical
value in the interval $(\scm,\scp)$.

\begin{definition}\label{def:flip-loci}
We define the \textit{flip loci} as
 \begin{align*}
 \cS_{\scp} &= \{ T\in\cN_{\scp} \ |
 \ \text{$T$ is $\scm$-unstable}\} \subset\cN_{\scp} \ ,\\
 \cS_{\scm} &= \{ T\in\cN_{\scm} \ |
 \ \text{$T$ is $\scp$-unstable}\}
 \subset\cN_{\scm} \ .
 \end{align*}
and $\cS_{\scpm}^s=\cS_\scpm \cap \cN_\scpm^s$ for the stable part of the
flip loci.
\end{definition}

Note that for $\s_c=\s_m$, $\cN_{\s_m^-}$ is empty, hence $\cN_{\smp}=
\cS_{\smp}$. Also $\cN_{\s_m}^s=\emptyset$, by the last part of
Proposition \ref{prop:sigma-coeff}. Analogously, when $n_1\neq n_2$,
$\cN_{\s_M^+}$ is empty, $\cN_{\s_M^-}= \cS_{\s_M^-}$ and $\cN_{\s_M}^s$
is empty.

\begin{lemma}\label{lem:fliploci}
 Let $\s_c$ be a critical value. Then
 \begin{itemize}
 \item[(1)] $\cN_{\scp}-\cS_{\scp}=\cN_{\scm}-\cS_{\scm}$.
 \item[(2)] $\cN^s_{\scp}-\cS_{\scp}^s=
   \cN^s_{\scm}-\cS_{\scm}^s=\cN^s_{\s_c}$.
 \end{itemize}
\end{lemma}

\begin{dem}
Item (1) is an easy consequence of the definition of flip loci. Item (2)
is the content of \cite[Lemma 5.3]{BGPG}.
\end{dem}

Now we want to describe the flip loci $\cS_{\s_c^{\pm}}$. Let $\s_c$ be a
critical value, and let $(n_1',n_2',d_1',d_2')$ such that $\lambda' \neq
\lambda$ and (\ref{eqn:sigmac}) holds. Put
$(n_1'',n_2'',d_1'',d_2'')=(n_1-n_1',n_2-n_2',d_1-d_1',d_2-d_2')$. Denote
$\cN_\s'=\cN_\s(n_1',n_2',d_1',d_2')$ and
$\cN_\s''=\cN_\s(n_1'',n_2'',d_1'',d_2'')$.

\begin{lemma} \label{lem:semistable}
Let $T\in \cS_\scp$ (resp.\ $T\in \cS_\scm$). Then $T$ is a
non-trivial extension
 \begin{equation}\label{eqn:extension-triples}
 0\to T'\to T\to T''\to 0,
 \end{equation}
where $\mu_{\s_c}(T')=\mu_{\s_c}(T)=\mu_{\s_c}(T'')$, $\lambda'<\lambda$
(resp.\ $\lambda'
 >\lambda$) and $T'$ and $T''$ are both
$\s_c$-semistable.

Conversely, suppose $T'\in \cN_{\s_c}'$ and $T''\in \cN_{\s_c}''$
are both $\s_c$-stable, and $\lambda'<\lambda$ (resp.\ $\lambda'
>\lambda$). Then for any non-trivial extension {\rm (\ref{eqn:extension-triples})}, $T$
lies in $\cS_\scp^s$ (resp.\ in $\cS_\scm^s$). Moreover, such $T$
can be written uniquely as an extension {\rm
(\ref{eqn:extension-triples})} with
$\mu_{\s_c}(T')=\mu_{\s_c}(T)$.

In particular, suppose $\s_c$ is not a critical value for the
moduli spaces of triples of types $(n_1',n_2',d_1',d_2')$ and
 $(n_1'',n_2'',d_1'',d_2'')$, $\GCD(n_1',n_2',d_1'+d_2')=1$ and
 $\GCD(n_1'',n_2'',d_1''+d_2'')=1$. Then if $\lambda'<\lambda$ (resp.\
 $\lambda'>\lambda$), there is a bijective correspondence
 between non-trivial extensions {\rm (\ref{eqn:extension-triples})},
with
 $T'\in \cN_{\s_c}'$ and $T''\in \cN_{\s_c}''$ and triples
 $T\in \cS_\scp$ (resp.\ $\Smenos$).
\end{lemma}

\begin{dem}
Let $T$ be a triple which is $\scp$-semistable but $\scm$-unstable. Then
Proposition \ref{prop:sigma-coeff} implies the existence of a subtriple
$T'$ and a quotient triple $T''=T/T'$ such that
$\mu_{\s_c}(T')=\mu_{\s_c}(T)=\mu_{\s_c}(T'')$ and $\lambda'<\lambda$. By
Remark \ref{rem:directsum}, $T$ defines a non-trivial extension in
$\Ext^1(T'',T')$. Now suppose $T'$ is $\s_c$-unstable; then there exists
$\tilde{T}\subset T'$ with
$\mu_{\s_c}(\tilde{T})>\mu_{\s_c}(T')=\mu_{\s_c}(T)$, contradicting the
$\s_c$-semistability of $T$. So $T'$ is $\s_c$-semistable. Also, if $T''$
is $\s_c$-unstable, then there exists a quotient $T''\surj \tilde{T}$ with
$\mu_{\s_c}(T'')=\mu_{\s_c}(T)>\mu_{\s_c}(\tilde{T})$, contradicting again
the $\s_c$-semistability of $T$. Hence $T''$ is $\s_c$-semistable.

The second statement is basically due to the uniqueness of Jordan-H\"older
filtrations. We provide the argument for completeness. Suppose
 $T$ is such a non-trivial extension and let us check that it is
 $\scp$-stable. For this, let
 $\tilde{T}\subset T$ be a subtriple. Consider the subtriple
 $\tilde{T}'$ given as the saturation of $\tilde{T}\cap T'$ (which
 is a triple of torsion-free subsheaves). Let
 $\tilde{T}''=\tilde{T}/\tilde{T}'$, so that there is an exact
 sequence
  $$
  \begin{array}{ccccccccc}
  0 & \to & \tilde{T}'& \to &\tilde{T}& \to &\tilde{T}''&\to
&0\, \, \\
    &    & \cap  &  &\cap  &  &\cap  \\
  0  & \to   & T' &  \to& T  &  \to&T'' &\to & 0\, .
 \end{array}
  $$
Clearly $\mu_{\s_c}(\tilde{T}')\leq \mu_{\s_c}(T')=\mu_{\s_c}(T)$
and $\mu_{\s_c}(\tilde{T}'')\leq \mu_{\s_c}(T'')=\mu_{\s_c}(T)$,
with a strict inequality if either of the subtriples
  \begin{equation} \label{eqn:trip-subtrip}
  \tilde{T}'\subset T',  \qquad
  \tilde{T}'' \subset T''
  \end{equation}
 is proper. In this case, Remark \ref{rem:slopes}
 gives that $\mu_{\s_c}(\tilde{T})<\mu_{\s_c}(T)$. By changing
 $\s_c$ to $\scp=\s_c+\varepsilon$, we still
 have $\mu_{\scp}(\tilde{T})<\mu_{\scp}(T)$.

 If both subtriples (\ref{eqn:trip-subtrip}) are not proper,
 then either $\tilde{T}'= T'$,
 $\tilde{T}'' =0 $ in which case $\tilde{T}=T'$ and hence
 $\mu_{\scp}(\tilde{T})<\mu_{\scp}(T)$ since $\lambda'<\lambda$;
 or $\tilde{T}'= 0$, $\tilde{T}'' =T''$ in which case there is a
 splitting $T''=\tilde{T}'' =\tilde{T}\inc T$ of the exact
 sequence $0\to T'\to T\to T''\to 0$, and the extension is
 trivial.

  The argument also shows that for such $T\in \cS_\scp$, the only
  subtriple with equal $\s_c$-slope is $T'$. This produces the
  uniqueness of the defining exact sequence.

  The last part follows since the conditions
  $\s_c$ not being a critical value for the moduli
spaces of triples of types $(n_1',n_2',d_1',d_2')$ and
 $(n_1'',n_2'',d_1'',d_2'')$, $\GCD(n_1',n_2',d_1'+d_2')=1$ and
 $\GCD(n_1'',n_2'',d_1''+d_2'')=1$ imply that $\cN'_{\s_c}$ and
 $\cN''_{\s_c}$ do not contain properly $\s_c$-semistable triples.
\end{dem}

\begin{theorem} \label{thm:Smas}
 Let $\s_c$ be a critical value with $\lambda'<\lambda$
 (resp.\ $\lambda'>\lambda$). Assume
 \begin{itemize}
 \item[(i)] $\s_c$ is not a critical value for the moduli spaces
 of triples of types $(n_1',n_2',d_1',d_2')$ and
 $(n_1'',n_2'',d_1'',d_2'')$, $\GCD(n_1',n_2',d_1'+d_2')=1$ and
 $\GCD(n_1'',n_2'',d_1''+d_2'')=1$.
 \item[(ii)] $\HH^0(C^\bullet(T'',T'))=\HH^2(C^\bullet(T'',T'))
  %\HH^0(C^\bullet(T',T''))=\HH^2(C^\bullet(T',T''))
 =0$, for every $(T',T'')\in \cN_{\s_c}'\times  \cN_{\s_c}''$.
 \item[(iii)] There are universal triples for $\cN_{\s_c}'$
 and $\cN_{\s_c}''$.
 \end{itemize}
 Then $\Smas$ (resp.\ $\Smenos$) is the projectivization of a
 bundle of
 rank $-\chi(T'',T')$ over $\cN_{\s_c}' \times \cN_{\s_c}''$.
\end{theorem}

\begin{dem}
Let $\cT'=(\cE_1',\cE_2',\Phi')\to \cN_{\s_c}'\times X$ and
$\cT''=(\cE_1'', \cE_2'',\Phi'')\to \cN_{\s_c}''\times X$ denote
universal triples, provided by Proposition
\ref{prop:triples-critical-range} (4). Let $B=\cN_{\s_c} '\x
\cN_{\s_c}''$ and pull back $\cT'$ and $\cT''$ to $B\x X$.
Consider the complex $C^\bullet(\cT'',\cT')$ as defined in
\eqref{eqn:extension-complex} and take relative hypercohomology
$\HH_\pi^i(C^\bullet(\cT'',\cT'))$ with respect to the projection
$\pi:B \x X\to B$. Setting
 $$
  W=\HH^1_\pi(C^\bullet(\cT'',\cT')),
 $$
we have an identification $\cS_\scp\iso \PP W$, by Lemma
\ref{lem:semistable} (that is, a bijection). By (ii), $W$ is a bundle of
rank $h^1(T'',T')=-\chi(T'',T')$.

We construct a universal triple for $\PP W$. Let $q:\PP W\to B$ be the
projection and consider the extension
 $$
 \begin{aligned}
 \xi \in \Ext^1(  (q \times 1_X)^* &\cT''  ,
 (q \times 1_X)^*\cT'\otimes \cO_{\PP
W}(1)) \\
&= \HH^1 (C^\bullet((q \times 1_X)^*\cT'',(q \times 1_X)^*\cT')
\otimes \cO_{\PP W}(1)) \\
 &= H^0 (\HH^1_\pi (C^\bullet(\cT'',\cT')) \otimes q_*\cO_{\PP
W}(1)) \\
 &= H^0(W \otimes W^*) =\End(W)
 \end{aligned}
 $$
corresponding to the identity homomorphism. (We have used that
$\HH^0_\pi (C^\bullet(\cT'',\cT'))=0$ in the third line of the
above chain of equalities.) This gives an extension
 $$
 0\to (q \times 1_X)^*\cT' \otimes \cO_{\PP W}(1) \to
 \cT \to (q \times 1_X)^*\cT'' \to 0,
 $$
which produces the required universal triple $\cT\to \PP W\times X$. This
defines an algebraic map $\PP W\to \cS_{\scp}$, which is the required
isomorphism.
\end{dem}

\begin{remark} \label{rem:parabolic}
This result is analogous to the results in \cite[Section 5]{GPGM} which
deal with the moduli space of parabolic triples. In that case, a suitable
choice of parabolic weights ensure that condition (i) in Theorem
\ref{thm:Smas} is always satisfied.

Condition (iii) in Theorem \ref{thm:Smas} is necessary for this
proof, and it is not automatic. For instance the moduli space of
stable bundles $M^s(n,d)$ does not have a universal bundle when
$\GCD(n,d)\neq 1$.
\end{remark}

The construction of the flip loci can be used for the critical
value $\s_c=\s_m=\mu_1-\mu_2$, which allows us to describe the
moduli space $\cN_{\smp}$. We refer to the value of $\s$ given by
$\s=\smp=\s_m+\varepsilon$ as \textit{small}.

\begin{proposition} \label{prop:moduli-small}
We have $\cN_{\s_m} \cong M(n_1,d_1)\times M(n_2,d_2)$ and only contains
strictly $\s_m$-semistable triples. There is a map
 $$
 \pi:\cN_{\s_m^+}(n_1,n_2,d_1,d_2) \to M(n_1,d_1) \times M(n_2,d_2)
 $$
which sends $T=(E_1,E_2,\phi)$ to $(E_1,E_2)$.
 \begin{enumerate}
 \item[(i)] If $\GCD(n_1,d_1)=1$, $\GCD(n_2,d_2)=1$ and
$\mu_1-\mu_2>2g-2$, then $\cN_{\s_m^+}^s=\cN_{\s_m^+}$ and it is
the projectivization of a bundle over $M(n_1,d_1) \times
M(n_2,d_2)$, of rank $n_2d_1-n_1d_2- n_1n_2(g-1)$.
 \item[(ii)] In general, if $\mu_1-\mu_2>2g-2$, the subset
 $$
 \pi^{-1}(M^s(n_1,d_1)\times
 M^s(n_2,d_2)) \subset \cN_{\s_m^+}
 $$
is a projective bundle over $M^s(n_1,d_1) \times M^s(n_2,d_2)$,
whose fibers are projective spaces of dimension $n_2d_1-n_1d_2-
n_1n_2(g-1)-1$.
 \end{enumerate}
\end{proposition}

\begin{dem}
Let $T=(E_1,E_2,\phi) \in \cN_{\s_m}$. Note that $0\to E_1$ is a
$\s_m$-destabilizing subtriple. As $T$ is polystable, $T=(E_1,0,0)
\oplus (0,E_2,0)$. (Here also $E_1$ and $E_2$ should be polystable
bundles.) Conversely if the triple $T=(E_1,E_2,\phi)$ has $\phi=0$
(and $n_1n_2\neq 0$), then automatically $T=(E_1,0,0)\oplus
(0,E_2,0)$ and it is $\s$-unstable for any $\s\neq \s_m$, by
Remark \ref{rem:directsum}. If $E_1$ and $E_2$ are polystable
bundles then $T\in \cN_{\s_m}$. This proves the first statement.

Next, the map $\pi$ is well-defined since a $\smp$-semistable
triple is automatically $\s_m$-semistable.

Now let us see item (i). Note that all the triples in
$\cN_{\s_m^+}$ are non-trivial extensions of $T'=(E_1,0,0)$ by
$T'' =(0,E_2,0)$ by Lemma \ref{lem:semistable}. The long exact
sequence in Proposition \ref{prop:hyper-equals-hom} yields
$\HH^0=0$, $\HH^1=\Hom(E_2,E_1)$ and $\HH^2=H^1(E_2^*\ox E_1)\cong
\Hom(E_1,E_2\otimes K)$, where we have abbreviated
$\HH^i=\HH^i(C^\bullet(T'',T'))$. Since both $E_1$ and $E_2$ are
stable bundles and $\mu_1>\mu_2+2g-2$, we get that $\HH^2=0$.
Therefore the conditions of Theorem \ref{thm:Smas} are satisfied
and $-\chi(T'',T')=-(n_1n_2(g-1) -n_2d_1+n_1d_2)$.

Item (ii) is similar, but now we do not have a universal bundle
for $M^s(n_1,d_1)$ or $M^s(n_2,d_2)$ at our disposal. Working in
the \'etale topology, we have (locally) a universal bundle which
yields that $\pi^{-1}(M^s(n_1,d_1) \times M^s(n_2,d_2))\to
M^s(n_1,d_1) \times M^s(n_2,d_2)$ is a fibration. The fiber over
$(E_1,E_2)$ is the projective space $\PP
\HH^1(C^\bullet(T'',T'))=\PP\Hom(E_2,E_1)$, using Lemma
\ref{lem:semistable}.
\end{dem}

%%%%%%%%%%%%%%%%%%%%%%%%%%%%%%%%%%%%%%%%%%%%%%%%%%%%%%%%%
\section{Moduli of triples of rank $(2,1)$ and pairs}
\label{sec:triples(2,1)}
%%%%%%%%%%%%%%%%%%%%%%%%%%%%%%%%%%%%%%%%%%%%%%%%%%%%%%%%%

Let $X$ be a smooth projective curve of genus $g\geq 2$. In this
section, we shall consider triples $T=(E_1,E_2,\phi)$ where $E_1$
is a vector bundle of degree $d_1$ and rank $2$ and $E_2$ is a
line bundle of degree $d_2$. Let $\moduli=\moduli(2,1,d_1,d_2)$
denote the moduli space of $\s$-polystable triples of such type.
By Proposition \ref{prop:alpha-range}, $\s$ is in the interval
  \begin{equation}\label{eqn:interv}
    I=[\s_m,\s_M]=
    [\mu_1-\mu_2\,,\,4(\mu_1-\mu_2)]=[d_1/2-d_2,2d_1-4d_2],
  \end{equation}
where $\mu_1-\mu_2\geq 0$. Otherwise $\cN_\s$ is empty.

\begin{theorem}\label{thm:moduli(2,1)}
For $\s\in I$, $\moduli$ is a projective variety. It is smooth and of
(complex) dimension $3g-2 + d_1 - 2 d_2$ at the stable points $\moduli^s$.
Moreover, for generic values of $\s$, $\moduli=\moduli^s$ (hence it is
smooth and projective).
\end{theorem}

\begin{dem}
The first line follows from Proposition \ref{prop:alpha-range}.
Now let $T=(E_1,E_2,\phi)$ be any $\s$-stable triple. Then
$\phi:E_2\to E_1$ is injective, since $E_2$ is a line bundle and
$\phi\neq 0$ ($\phi=0$ would imply that $T$ is decomposable, hence
not $\s$-stable). Therefore by Theorem \ref{thm:smoothdim} (4),
$\moduli$ is smooth at $T$ and of the stated dimension. Finally,
since $\GCD(1,2,d_1+d_2)=1$, Proposition
\ref{prop:triples-critical-range} (3) implies the last assertion.
\end{dem}

The moduli spaces of $\s$-stable triples of rank $(2,1)$ are
intimately related to the moduli space of $\tau$-stable pairs of
rank $2$. A holomorphic pair $(E,\phi)$ is formed by a rank $n$
holomorphic bundle $E$ and a non-zero holomorphic section $\phi\in
H^0(E)$. There is a notion of stability depending on a real
parameter $\tau\in\RR$ given in \cite[Definition 4.7]{GP}. A
holomorphic pair $(E,\phi)$ is $\tau$-stable if

\begin{enumerate}
 \item[(1)] $\mu(E')<\t$ for every subbundle $E' \subset E$ with
 $\rk(E')>0$.
  \item[(2)] $\mu(E/E')>\t$ for every subbundle $E' \subset E$ with
 $0<\rk(E')<r$ and $\phi\in H^0(E')$.
 \end{enumerate}

There is a moduli space $\frM_\tau(n,d)$ of $\tau$-stable pairs of
rank $n$ and degree $d$, and a moduli space $\frM_\tau(n,\Lambda)$
of $\tau$-stable pairs of rank $n$ and satisfying that the
determinant of the bundle $E$ is some fixed line bundle $\Lambda$
of degree $d$ on $X$. These moduli spaces have been studied in
\cite{B, BD, GP}.

Bertram \cite{B} and Thaddeus \cite{Th} gave a explicit GIT
construction of the moduli spaces $\frM_\tau(2,\Lambda)$ of pairs
$(E,\phi)$ of rank $2$ with fixed determinant $\det (E)=\Lambda$,
where $\Lambda$ is a fixed line bundle.

We have the following relationship between the moduli spaces of
triples of rank $(2,1)$ and the moduli spaces of pairs.

\begin{lemma}\label{lem:triplesversuspairs}
 We have an isomorphism $\cN_\s(2,1,d_1,d_2)\cong \frM_\tau(2,d_1-2d_2)
\times \Jac^{d_2}(X)$.
\end{lemma}

\begin{dem}
In fact, a triple $(E_1,E_2,\phi)$ is $\s$-stable if and only if
the pair $(E_1\otimes E_2^*,\phi)$ is $\tau$-stable with $\tau=
\frac13 (\s+d_1-2 d_2)$ (cf.\ \cite[Definition 4.27]{GP} and
\cite[Proposition 3.4]{BGP}). Therefore there exists a bijective
correspondence $\cN_\s(2,1,d_1,d_2)\cong \frM_\tau(2,d_1-2d_2)
\times \Jac^{d_2}(X)$, given by $(E_1,E_2,\phi) \mapsto
((E_1\otimes E_2^*,\phi),E_2)$. This is clearly an isomorphism.
\end{dem}

\subsection*{Critical values and flip loci}
Let us compute the critical values corresponding to $n_1=2$, $n_2=1$.
Following Definition \ref{def:critical}, we have the following
possibilities:

\begin{enumerate}
 \item[(1)] $n'_1=1$, $n'_2=0$. The corresponding $\s_c$-destabilizing
  subtriple is of the form $0\to E_1'$, where $E_1'=M$ is a line
bundle of
  degree $\deg(M)=d_M$. The critical value is
  $$
   \s_c= 3d_M-d_1-d_2\, .
  $$
 \item[(2)] $n'_1=1$, $n'_2=1$. The corresponding
 $\s_c$-destabilizing subtriple $T'$ is of the form $E_2\to
 E_1'$, where $E_1'$ is a line bundle. Let $T''=T/T'$ be the
 quotient bundle, which is of the form $0\to E_1''$, where
 $E_1''=M$ is a line bundle, and let $d_M=\deg(M)$ be its degree.
 Then $d_2'=d_2$, $d_1'=d_1-d_M$ and
  $$
   \s_c=-\big( 3(d_1-d_M+d_2)-2(d_1+d_2)\big)= 3d_M-d_1-d_2\, .
  $$

 \item[(3)] $n'_1=2$, $n'_2=0$. In this case, the only possible
subtriple
  is $0\to E_1$. This produces the critical value
  $$
   \s_c=\frac{d_1-2d_2}{2}= \mu_1-\mu_2=\s_m\, ,
  $$
  i.e.,\ the minimum of the interval $I$ for
  $\s$.
 \item[(4)] $n'_1=0$, $n'_2=1$. The subtriple $T'$ must be of the
form $E_2\to
  0$. This forces $\phi=0$ in $T=(E_1,E_2,\phi)$. So $T$ is
decomposable,
  of the form $T'\oplus T''=(0,E_2,0)\oplus (E_1,0,0)$. By Remark
  \ref{rem:directsum}, $T$ is $\s$-unstable for any
  $\s\neq \s_c$, where
  $$
   \s_c=\frac{2d_2-d_1}{-2}= \mu_1-\mu_2=\s_m\, .
  $$
  Actually $T$ is $\s_m$-semistable if and only if $E_1$ is a rank $2$
  semistable bundle.
\end{enumerate}

\begin{lemma} \label{lem:dM}
Let $\s_c=3d_M-d_1-d_2$ be a critical value. Then
 \begin{equation*}
   \mu_1\leq d_M\leq d_1-d_2 \, ,
 \end{equation*}
and $\s_c=\s_m \Leftrightarrow d_M=\mu_1$.
\end{lemma}

\begin{dem}
This simply consists of rewriting the inequalities $\s_m\leq
\s_c\leq \s_M$.
\end{dem}

%%%%%%%%%%%%%%%%%%%%%%%%%%%%%%%%%%%%%%%%%%%%%%%%%%%%%%%%%
\section{Hodge polynomials of the moduli spaces of
triples of rank $(2,1)$} \label{sec:Hodge(2,1)}
%%%%%%%%%%%%%%%%%%%%%%%%%%%%%%%%%%%%%%%%%%%%%%%%%%%%%%%%%

In this section we are going to compute the Hodge polynomial of
the moduli space of $\s$-stable triples of type $(2,1,d_1,d_2)$
without making use of the blow-ups and blow-downs constructed by
Thaddeus in \cite{Th}. As $\cN_\s=\cN_\s(2,1,d_1,d_2)$ is smooth
and projective for non-critical values $\s$, by Theorem
\ref{thm:moduli(2,1)}, the polynomial $e(\cN_\s)$ is the usual
Hodge polynomial of $\cN_\s$.

Let $\s_c$ be a critical value. By Lemma \ref{lem:fliploci}, we
have that
 $$
  \cN_{\scp} - \cS_{\scp}= \cN_{\scm}- \cS_{\scm}\, .
 $$
Now Theorem \ref{thm:Du} implies that the Hodge polynomials
satisfy
 \begin{equation}\label{eqn:acabando}
  e(\cN_{\scp})-e(\cS_{\scp})=
  e( \cN_{\scm})-e(\cS_{\scm})\, .
 \end{equation}
We have the following computation for the difference
$e(\cS_{\scm})- e( \cS_{\scp})$.

\begin{lemma}\label{lem:flip-contribution}
 Let $\s_c=3d_M-d_1-d_2$ be a critical value for $\moduli(2,1,d_1,d_2)$ and assume that $\s_c\neq
\s_m$.
 Then $e(\cS_\scm)-e(\cS_\scp)$ is
  $$
  \coeff_{x^0} \left[
    \frac{((uv)^{d_1-d_2-d_M}-(uv)^{2d_M-d_1+g-1})
    (1+u)^{2g}(1+v)^{2g}
    (1+ux)^{g}(1+vx)^{g}}
    {(1-uv)(1-x)(1-uvx)x^{d_1-d_2-d_M}}
  \right]
  \, .
  $$
%
%$$
%\frac{((uv)^{2d_1-2d_2-2d_M}-(uv)^{4d_M-2d_1+2g-2})
%      (1+u)^{2g}(1+v)^{2g}}{(1-uv)}
%  \coeff_{x^0} \left[
%    \frac{(1+ux)^{g}(1+vx)^{g}}
%    {(1-x)(1-uvx)x^{d_1-d_2-d_M}}
%  \right]\,.
%$$
\end{lemma}

\begin{dem}
Since $\s_c\neq \s_m$, the critical value $\s_c=3d_M-d_1-d_2$ determines
the value of $d_M$ and there are only two possibilities for flip loci:
 \begin{enumerate}
 \item[(1)] The subtriple $T'$ is of the form $0\to E_1'$ with
$E_1'=M$ a line bundle
  of degree $d_M$, and the quotient triple $T''$ is of the form
$E_2\to E_1''$
  where $E_1''$ is a line bundle of degree $d_1-d_M$. Then by
  Lemmas \ref{lem:(n,0)} and \ref{lem:type(1,1)},
   \begin{align*}
   \cN'_{\s_c} &= \cN_{\s_c}(1,0,d_M,0)= \Jac^{d_M}X \, ,\\
   \cN''_{\s_c} &= \cN_{\s_c}(1,1,d_1-d_M,d_2) =\Jac^{d_2}X \x
   \Sym^{d_1-d_M-d_2}X\, .
   \end{align*}
  Here note that the second line needs that $\s_c\neq
  d_1-d_M-d_2$ (which is the only critical value for the moduli spaces
  of triples of rank $(1,1)$. This is true since $\mu_1<d_M$).
  Now note that $\lambda'=0<\lambda=\frac13$,
  $\HH^0(C^\bullet(T'',T'))=0$ by Proposition \ref{prop:h0-vanishing}.
  As $\phi'':E_2''\to E_1''$ is injective, $E_1^{''*}\otimes
  E_1'\to E_2^{''*}\otimes E_1'$ is generically surjective, so
$H^1(E_1^{''*}\otimes
  E_1')\to H^1(E_2^{''*}\otimes E_1')$ is surjective. Hence
  the long exact sequence in
  Proposition \ref{prop:hyper-equals-hom} implies that
  $\HH^2(C^\bullet(T'',T'))=0$.
  Therefore, by Theorem \ref{thm:Smas}, $\cS_\scp$ is the
  projectivization of  a rank $-\chi(T'',T')$ bundle over
  $\cN'_{\s_c}\times \cN''_{\s_c}$, where
  $$
    -\chi(T'',T')=d_1-d_M-d_2,
  $$
  using Proposition \ref{prop:chi(T'',T')}. Now Lemma \ref{lem:vb}
  yields
  \begin{align*}
  e(\cS_\scp)&=
  e(\PP^{d_1-d_M-d_2-1}) e(\cN'_{\s_c}\times
  \cN''_{\s_c}) \\
  &=e_{d_1-d_M-d_2}\,
  e(\Jac\ X)^2 e(\Sym^{d_1-d_M-d_2}X)\, ,
  \end{align*}
  with the notation in (\ref{eqn:Pn}). This formula also holds
  when $d_1-d_M-d_2=0$, since in such case $\cS_{\scp}=\emptyset$
  and $e_{d_1-d_M-d_2}=0$ (note that it cannot happen that
  $d_1-d_M-d_2<0$).

\item[(2)] The subtriple $T'$ is of the form $E_2\to E_1'$ with
quotient of the form $0\to E_1''$, where $E_1''=M$ is a line bundle of
degree $d_M$. Now
\begin{align*}
    \cN'_{\s_c} &= \cN_{\s_c}(1,1,d_1-d_M,d_2) =\Jac^{d_2}X \x
    \Sym^{d_1-d_M-d_2}X\, , \\
    \cN''_{\s_c} &= \cN_{\s_c}(1,0,d_M,0)= \Jac^{d_M}X \, .
\end{align*}
Again $\HH^0(C^\bullet(T'',T'))=0$ by Proposition
\ref{prop:h0-vanishing}, and $\HH^2(C^\bullet(T'',T'))=0$ by the
long exact sequence in Proposition \ref{prop:hyper-equals-hom} and
using $E_2''=0$. As $\lambda'=\frac12>\lambda=\frac13$, Theorem
\ref{thm:Smas} says that $\cS_\scm$ is the projectivization of  a
rank $-\chi(T'',T')$ bundle over $\cN'_{\s_c}\times \cN''_{\s_c}$,
where
$$
  -\chi(T'',T')=2d_M-d_1+g-1\, ,
$$
by Proposition \ref{prop:chi(T'',T')}. So
\begin{equation*}
  e(\cS_\scm)=e_{2d_M-d_1+g-1}\, e(\Jac\ X)^2 e(\Sym^{d_1-d_M-d_2}X)\, .
\end{equation*}
\end{enumerate}

Using (\ref{eqn:Jac}) for the Hodge polynomial of the Jacobian of
$X$ and (\ref{eqn:sym}) for the Hodge polynomial of the symmetric
product of $X$, we get
 $$
 \begin{aligned}
    e(\cS_\scm)-e(\cS_\scp) = &( e_{2d_M-d_1+g-1}-e_{d_1-d_M-d_2})
    e(\Jac\ X)^2 e(\Sym^{d_1-d_M-d_2}X) \\
    = &
    \left(
    \frac{(uv)^{d_1-d_2-d_M}-(uv)^{2d_M-d_1+g-1}}{1-uv}\right)
    (1+u)^{2g}(1+v)^{2g}\\
     & \qquad \quad \cdot \,
   \coeff_{x^0}\frac{(1+ux)^{g}(1+vx)^{g}}{(1-x)(1-uvx)x^{d_1-d_2-d_M}}\,
     ,
 \end{aligned}
 $$
from which the stated result is obtained.
\end{dem}

\begin{theorem} \label{thm:polinomio(2,1)no-critico}
Let $X$ be a smooth projective curve of genus $g \ge 2$  and
consider the moduli space $\moduli=\moduli (2,1,d_1,d_2)$, for a
non-critical value $\s>\s_m$. Set
$d_0=\Big[\frac13(\s+d_1+d_2)\Big]+1$. Then the Hodge polynomial
of $\moduli$ is
 $$
  e(\moduli)= \coeff_{x^0}
\left[\frac{(1+u)^{2g}(1+v)^{2g}(1+ux)^{g}(1+vx)^{g}}{(1-uv)(1-x)(1-uvx)x^{d_1-d_2-d_0}}
    \Bigg(\frac{(uv)^{d_1-d_2-d_0}}{1-(uv)^{-1}x}-
    \frac{(uv)^{-d_1+g-1+2d_0}}{1-(uv)^2x}\Bigg)\right] .
 $$
\end{theorem}

\begin{dem}
By (\ref{eqn:acabando}), we have that $e(\cN_{\scm})- e(
\cN_{\scp})=e(\cS_{\scm})- e( \cS_{\scp})$. Lemma
\ref{lem:flip-contribution} implies that this quantity equals
 $$
 \coeff_{x^0} \left[
    \frac{((uv)^{d_1-d_2-d_M}-(uv)^{2d_M-d_1+g-1})
 (1+u)^{2g}(1+v)^{2g}(1+ux)^{g}(1+vx)^{g}}{(1-uv)(1-x)(1-uvx)x^{d_1-d_2-d_M}}
  \right]
  \, .
 $$
Now we add up all the contributions for $\s<\s_c\leq \s_M$. This
corresponds to
 \begin{equation}\label{eqn:star}
 \frac13(\s+d_1+d_2)<d_M\leq d_1-d_2.
 \end{equation}
Note that since $\s$ is not a critical value, we cannot have
equality in the left, hence the left hand term is not an integer.
Therefore (\ref{eqn:star}) is equivalent to $d_0\leq d_M\leq
d_1-d_2$, with $d_0$ given as in the statement. Thus
 $$
 \begin{aligned}
   e(\cN_\s) &=\coeff_{x^0}
\Bigg(\frac{(1+u)^{2g}(1+v)^{2g}(1+ux)^{g}(1+vx)^{g}(uv)^{d_1-d_2}}{(1-uv)(1-x)(1-uvx)x^{d_1-d_2}}
    \sum_{d_M=d_0}^{d_1-d_2} (uv)^{-d_M}x^{d_M}\\
    &\qquad\qquad\quad
    -\,
\frac{(1+u)^{2g}(1+v)^{2g}(1+ux)^{g}(1+vx)^{g}(uv)^{-d_1+g-1}}{(1-uv)(1-x)(1-uvx)x^{d_1-d_2}}
    \sum_{d_M=d_0}^{d_1-d_2} (uv)^{2d_M}x^{d_M}
    \Bigg)\, .
 \end{aligned}
 $$
If we add to the summations terms with $d_M>d_1-d_2$, then all the
new terms contribute positive powers of $x$ to the global
expression, so they do not appear after extracting the coefficient
of $x^0$. This means that we can take the sum from $d_0$ to
$\infty$. Using
 $$
    \sum_{d_M=d_0}^\infty (uv)^{-d_M}x^{d_M}
    =\frac{(uv)^{-d_0}x^{d_0}}{1-(uv)^{-1}x}
    \quad\mbox{ and }\quad
    \sum_{d_M=d_0}^\infty (uv)^{2d_M}x^{d_M}
    =\frac{(uv)^{2d_0}x^{d_0}}{1-(uv)^{2}x}\, ,
 $$
we get
 $$
    \begin{aligned}
    e(\moduli)&=
    \coeff_{x^0}
\Bigg(\frac{(1+u)^{2g}(1+v)^{2g}(1+ux)^{g}(1+vx)^{g}(uv)^{d_1-d_2}(uv)^{-d_0}
    x^{d_0}}
    {(1-uv)(1-x)(1-uvx)(1-(uv)^{-1}x)x^{d_1-d_2}}\\
    &\qquad\qquad\quad
    -\,
\frac{(1+u)^{2g}(1+v)^{2g}(1+ux)^{g}(1+vx)^{g}(uv)^{-d_1+g-1}(uv)^{2d_0}
    x^{d_0}}
    {(1-uv)(1-x)(1-uvx)(1-(uv)^{2}x)x^{d_1-d_2}}
    \Bigg) \, ,
    \end{aligned}
 $$
from which the result follows.
\end{dem}

Obviously, $e(\moduli)=0$ for $\s<\s_m$, since in this case
$\moduli$ is empty.

We recover the Poincar{\'e} polynomial of the moduli space
$\cN_\s(2,1,d_1,d_2)$. This agrees with the formula given in
\cite[Theorem 6.4]{GPGM} for the case of parabolic triples, if we
consider in the latter that there are no parabolic points.

\begin{corollary} \label{cor:polinomio(2,1)critico}
Let $X$ be a smooth projective curve of genus $g \ge 2$  and
consider the moduli space $\moduli=\moduli (2,1,d_1,d_2)$ for a
non-critical value $\s>\s_m$. Set
$d_0=\Big[\frac13(\s+d_1+d_2)\Big]+1$. Then the Poincar{\'e}
polynomial of $\moduli=\moduli (2,1,d_1,d_2)$ is
 $$
  P_t(\moduli)=
 \coeff_{x^0}
 \left[\frac{(1+t)^{4g}(1+tx)^{2g}}{(1-t^2)(1-x)(1-t^2x)x^{d_1-d_2-d_0}}
    \Bigg(\frac{t^{2d_1-2d_2-2d_0}}{1-t^{-2}x}-
    \frac{t^{-2d_1+2g-2+4d_0}}{1-t^4x}\Bigg)\right] .
 $$
\end{corollary}

\begin{dem}
 Substitute $u=t$, $v=t$ into the formula of Theorem
 \ref{thm:polinomio(2,1)no-critico}.
\end{dem}

\medskip

\noindent\textbf{Proof of Theorem \ref{thm:main}:}\/ By Lemma
\ref{lem:triplesversuspairs}, we have an isomorphism
 $$
    \frM_\tau(2,d) \times
    \Jac^{d_2}(X) \cong \cN_\s(2,1,d_1,d_2),
 $$
with $d_1=d+2d_2$ and $\s=3\t-d$. Then
 $$
  e(\frM_\t(2,d))= \frac{e(\cN_\s(2,1,d_1,d_2))}{(1+u)^g(1+v)^g}\,
  .
 $$
So Theorem \ref{thm:polinomio(2,1)no-critico} gives
 \begin{eqnarray*}
  e(\frM_{\t}(2,d)) &=& \coeff_{x^0}
  \left[\frac{(1+u)^{g}(1+v)^{g}(1+ux)^{g}(1+vx)^{g}}{(1-uv)(1-x)(1-uvx)x^{d+d_2-d_0}}
    \Bigg(\frac{(uv)^{d+d_2-d_0}}{1-(uv)^{-1}x}-
    \frac{(uv)^{g-1-d-2d_2+2d_0}}{1-(uv)^2x}\Bigg)\right] \\
    &=& \coeff_{x^0}
  \left[\frac{(1+u)^{g}(1+v)^{g}(1+ux)^{g}(1+vx)^{g}}{(1-uv)(1-x)(1-uvx)x^{d-1-[\t]}}
    \Bigg(\frac{(uv)^{d-1-[\t]}}{1-(uv)^{-1}x}-
    \frac{(uv)^{g+1-d+2[\t]}}{1-(uv)^2x}\Bigg)\right],
 \end{eqnarray*}
since $d_0=\left[\frac13 (\s+d_1+d_2)\right]+1=[\t]+d_2+1$. Note
that $\t$ is not an integer and
  $$
  \t\in J=[d/2,d],
  $$
using (\ref{eqn:interv}), $d_1=d+2d_2$ and $\t=\frac13(\s+d)$. The
critical values for $\frM_\t(2,d)$ are the integers in $J$.

For the second part, note that the determinant map gives a
fibration
 \begin{equation}\label{eqn:fibr}
 \frM_\t(2,d) \to \Jac^d X
 \end{equation}
with fibers isomorphic to $\frM_\t(2,\Lambda)$. Note that all the
moduli spaces $\frM_\t(2,\Lambda')$ for any $\Lambda'\in \Jac^d X$
are isomorphic to each other and smooth (for non-critical value
$\tau$). The Serre spectral sequence of the fibration has $E_2$
term isomorphic to $H^*(\frM_{\t}(2,\Lambda))\otimes H^*(\Jac^d
X)$ and converges to $H^*(\frM_{\t}(2,d))$. By \cite[formula
(4.1)]{Th} (upon the substitution $i=d-1-[\tau]$),
  $$
  P_t(\frM_{\t}(2,\Lambda)) = \coeff_{x^0}
  \left[\frac{(1+tx)^{2g}}{(1-t^2)(1-x)(1-t^2x)x^{d-1-[\t]}}
    \Bigg(\frac{t^{2d-2-2[\t]}}{1-t^{-2}x}-
    \frac{t^{2g+2-2d+4[\t]}}{1-t^4x}\Bigg)\right],
  $$
and by our above calculation
  $$
  P_t(\frM_{\t}(2,d)) = \coeff_{x^0}
  \left[\frac{(1+t)^{2g}(1+tx)^{2g}}{(1-t^2)(1-x)(1-t^2x)x^{d-1-[\t]}}
    \Bigg(\frac{t^{2d-2-2[\t]}}{1-t^{-2}x}-
    \frac{t^{2g+2-2d+4[\t]}}{1-t^4x}\Bigg)\right],
  $$
from where we obtain $P_t(\frM_{\t}(2,d))=P_t(\frM_{\t}(2,\Lambda)
)P_t(\Jac X)$. Therefore the Serre spectral sequence degenerates,
and so the fibration (\ref{eqn:fibr}) is a rationally
cohomologically trivial fibration. This implies that there is an
isomorphism of Hodge structures
 $$
 H^*(\frM_{\t}(2,d)) \cong H^*(\frM_{\t}(2,\Lambda))\otimes H^*(\Jac^d X) .
 $$
So
 $$
 e(\frM_{\t}(2,d))=e(\Jac^d X) e(\frM_{\t}(2,\Lambda)).
 $$
The result follows from this. \hfill $\Box$

%%%%%%%%%%%%%%%%%%%%%%%%%%%%%%%%%%%%%%%%%%%%%%%%%%%%%%%%%%%%%%
\section{Hodge polynomial of the moduli space of triples of
rank $(1,2)$} \label{sec:Hodge(1,2)}
%%%%%%%%%%%%%%%%%%%%%%%%%%%%%%%%%%%%%%%%%%%%%%%%%%%%%%%%%%%%%%

For completeness, we include the computation of the Hodge
polynomial of the moduli of triples of rank $(1,2)$. Such triples
are of the form $\phi: E_2\to E_1$, where $E_2$ is a rank $2$
bundle and $E_1$ is a line bundle. By Proposition
\ref{prop:alpha-range}, $\s$ is in the interval
  $$
    I=[\s_m,\s_M]=
    [\mu_1-\mu_2\,,\,4(\mu_1-\mu_2)]=[d_1-d_2/2,4d_1-2d_2],
  $$
where $\mu_1-\mu_2\ge 0$. Otherwise $\cN_\s$ is empty. By the
duality result of Proposition \ref{prop:duality}, one has an
isomorphism
 $$
 \moduli(1,2,d_1,d_2) \cong \moduli(2,1,-d_2,-d_1).
 $$
So the study of these triples reduces to the case of triples of
rank $(2,1)$. From Theorem \ref{thm:moduli(2,1)} one immediately
obtains

\begin{theorem}\label{thm:moduli(1,2)}
For $\s\in I$, $\moduli$ is a projective variety. It is smooth and of
(complex) dimension $3g -2 + 2 d_1 - d_2$ at the stable points
$\moduli^s$. Moreover, for generic values of $\s$, $\moduli=\moduli^s$
(hence it is smooth and projective). %\hfill $\Box$
\end{theorem}

{}From Lemma \ref{lem:dM}, we obtain

\begin{lemma} \label{lem:dM-2}
The critical values for $\moduli(1,2,d_1,d_2)$ are the numbers
$\s_c=3d_M+d_1+d_2$, where $-\mu_2\leq d_M \leq d_1-d_2$. Also $\s_c=\s_m
\Leftrightarrow d_M=-\mu_2$.
\end{lemma}

\begin{theorem} \label{thm:-polinomiono(1,2)no-critico}
Consider $\cN_\s=\cN_\s(1,2,d_1,d_2)$. Let $\s>\s_m$ be a non-critical
value. Set $d_0=\Big[\frac13(\s-d_1-d_2)\Big]+1$. Then the Hodge
polynomial of $\moduli$ is
 $$
  e(\moduli)=
  \coeff_{x^0}
  \left[\frac{(1+u)^{2g}(1+v)^{2g}(1+ux)^{g}(1+vx)^{g}}
            {(1-uv)(1-x)(1-uvx)x^{d_1-d_2-d_0}}
    \Bigg(\frac{(uv)^{d_1-d_2-d_0}}{1-(uv)^{-1}x}-
    \frac{(uv)^{d_2+g-1+2d_0}}{1-(uv)^2x}\Bigg)\right] .
 $$
\end{theorem}

\begin{dem}
We use that $e(\moduli(1,2,d_1,d_2))=e(\moduli(2,1,-d_2,-d_1))$
and the formula in Theorem \ref{thm:polinomio(2,1)no-critico},
where $d_1$ and $d_2$ are substituted by $-d_2$, $-d_1$ and
 $$
 d_0= \left[ \frac13 (\s -d_2-d_1)\right] +1 \, .
 $$
\mbox{}\vskip-3\baselineskip\mbox{}\hfill\end{dem}

%%%%%%%%%%%%%%%%%%%%%%%%%%%%%%%%%%%%%%%%%%%%%%%%%%%%%%%%%%%%%%
\section{Hodge polynomial of the moduli of bundles} \label{sec:poly-bundles-rk-two}
%%%%%%%%%%%%%%%%%%%%%%%%%%%%%%%%%%%%%%%%%%%%%%%%%%%%%%%%%%%%%%

Let $M(2,d)$ denote the moduli space of polystable vector bundles
of rank $2$ and degree $d$ over $X$. Note that $M(2,d)\cong
M(2,d+2k)$, for any integer $k$, so there are basically two moduli
spaces, depending on whether the degree is even or odd. We are
going to apply the results of the Section \ref{sec:Hodge(2,1)} to
compute the Hodge polynomials of $M(2,d)$ and $M(2,\Lambda)$ with
$d$ odd, and $\Lambda$ a fixed line bundle of degree $d$ on $X$,
recovering the results of \cite{EK} for $M(n,d)$ and for
$M(n,\Lambda)$, when $n=2$. Note that the Hodge polynomial of the
moduli space $M(2,\Lambda)$ was also found by del Ba\~{n}o
\cite[$\S 3$]{BaR} using motivic techniques. The formula for the
Hodge polynomial of $M(2,\Lambda)$ can also be obtained from the
previous paper \cite{Ne1}, where an explicit basis for the
cohomology groups of $M(2,\Lambda)$ is given and the elements of
that basis are of pure Hodge type provided one chooses the
underlying basis of $H^1(X)$ to consist of elements of pure Hodge
type.

\begin{proposition}\label{prop:rank2odd}
The Hodge polynomial of $M(2,d)$, for odd degree $d$, is
 $$
    e(M(2,d)) = \frac{(1+u)^{g}(1+v)^{g}
                    (1+u^2v)^{g}(1+uv^2)^{g}
                    -(uv)^{g}(1+u)^{2g}(1+v)^{2g}}
                {(1-uv)(1-(uv)^2)}\,.
 $$
Suppose that $\Lambda$ is a fixed line bundle of odd degree on
$X$. Then the Hodge polynomial of $M(2,\Lambda)$ is
 $$
    e(M(2,\Lambda)) = \frac{
                    (1+u^2v)^{g}(1+uv^2)^{g}
                    -(uv)^{g}(1+u)^{g}(1+v)^{g}}
                {(1-uv)(1-(uv)^2)}\,.
 $$
\end{proposition}

\begin{dem}
Atiyah and Bott \cite[Proposition 9.7]{AB} show that the
determinant map $M(n,d) \to \Jac^d X$ with fibre $M(n,\Lambda)$ is
a cohomologically trivial fibration (for the case $n=2$, the
cohomological triviality of this fibration was proved in
\cite{Ne2}). This implies \cite[Lemma 3]{EK} an isomorphism of
Hodge structures
 $$
 H^*(M(n,d)) \cong H^*(M(n,\Lambda)) \otimes H^*(\Jac^d X)\, ,
 $$
from where, in particular,
 \begin{equation*}\label{eqn:isomorphismohodge}
 e(M(2,d))=e(M(2,\Lambda))e(\Jac\,X)\, .
 \end{equation*}
Therefore the second formula is a consequence of the first. We
shall prove the first one.

Choose $(n_1,d_1)=(2,d)$ and $(n_2,d_2)=(1,d_2)$. If $d_2$ is very
negative then $\mu_1-\mu_2=d/2 - d_2>2g-2$. We shall choose the
minimum possible value $d-2d_2=4g-3$. Then Proposition
\ref{prop:moduli-small} (i) applies and $\cN_{\smp}$ is the
projectivization of a vector bundle over $M(2,d)\x \Jac^{d_2} X$,
of rank $d-2d_2-2(g-1)=2g-1$. Lemma \ref{lem:vb} implies that
 \begin{equation}\label{eqn:solve-for-M(2,d)}
 e(\modulimmas)= e(\Jac^{d_2} X)\, e(M(2,d))\, e_{2g-1}\, ,
 \end{equation}
where
 \begin{equation}\label{eqn:solve-for-M(2,d)2}
 e(\Jac^{d_2} X)=(1+u)^{g}(1+v)^{g}\, \quad \text{and}
 \quad e_{2g-1}= \frac{1-(uv)^{2g-1}}{1-uv}\, .
 \end{equation}
To compute the left hand side of (\ref{eqn:solve-for-M(2,d)}), we use
Theorem \ref{thm:polinomio(2,1)no-critico} for
$\s=\smp=\mu_1-\mu_2+\varepsilon$, $\varepsilon>0$ small. Clearly,
 $$
 d_0=\big[\mbox{$\frac13$}(\mu_1-\mu_2+\varepsilon+2\mu_1+\mu_2)\big]+1=
 [\mu_1]+1=\frac{d+1}2 \,.
 $$
The Hodge polynomial of $\modulimmas$ is thus
 $$
    \begin{aligned}
    e(\modulimmas)&=
    \coeff_{x^0}
  \Bigg[\frac{(1+u)^{2g}(1+v)^{2g}(1+ux)^{g}(1+vx)^{g}}
    {(1-uv)(1-x)(1-uvx)x^{d-d_2-(d+1)/2}} \Bigg(
   \frac{(uv)^{d-d_2-d_0}}{1-(uv)^{-1}x}
   - \frac{(uv)^{-d+g-1+d_0}}{1-(uv)^{2}x} \Bigg) \Bigg] \, .
    \end{aligned}
 $$
Plugging this into (\ref{eqn:solve-for-M(2,d)}) and using
(\ref{eqn:solve-for-M(2,d)2}), we get
 $$
 \begin{aligned}
    e(M(2,d))= \frac{(1+u)^{g}(1+v)^{g}}{1-(uv)^{2g-1}}
    \coeff_{x^0}\Bigg[ &
    \frac{(1+ux)^{g}(1+vx)^{g}}{(1-x)(1-uvx)x^{2g-2}} \Bigg(
    \frac{(uv)^{2g-2}}{1-(uv)^{-1}x} - \frac{(uv)^{g}}{1-(uv)^2x} \Bigg) \Bigg]\, .
 \end{aligned}
 $$

To compute the coefficient of $x^0$, we work as follows. Denote
$$
\begin{aligned}
    f(x)&=(1+ux)^{g}(1+vx)^{g}x^{1-2g}\, , \\
    F(a,b,c)&=\coeff_{x^{0}}
    \frac{xf(x)}{(1-ax)(1-bx)(1-cx)}\, ,
\end{aligned}
$$
with $a,b,c$ distinct and different from zero. Then
\begin{equation}\label{Hodge-and-residue}
    e(M(2,d))=\frac{(1+u)^{g}(1+v)^{g}}{1-(uv)^{2g-1}}
        \bigg[
        F(1,uv,{\textstyle\frac1{uv}})\,(uv)^{2g-2}
        -F(1,uv,(uv)^2)\,(uv)^g
        \bigg]\,.
\end{equation}
On the other hand, we have the equality
 $$
    \frac{1}{(1-ax)(1-bx)(1-cx)}=
    \frac{A}{1-ax}+\frac{B}{1-bx}+\frac{C}{1-cx}\,,
 $$
where
 $$
  A =\frac{a^2}{(a-b)(a-c)}\,,\quad
    B=\frac{b^2}{(b-a)(b-c)}\,,\quad
    C=\frac{c^2}{(c-a)(c-b)}\,,
 $$
and so, by the residue theorem applied to
$\frac{f(x)}{(1-ax)(1-bx)(1-cx)}$\ ,
 \begin{equation}\label{eqn:F}
    F(a,b,c)
    =-\Res\left\{\frac{A\,f(x)}{1-ax};\frac1a\right\}
    -\Res\left\{\frac{B\,f(x)}{1-bx};\frac1b\right\}
    -\Res\left\{\frac{C\,f(x)}{1-cx};\frac1c\right\}\, ,
 \end{equation}
because $f(x)$ is holomorphic in $\CC-\{0\}$ and
$\frac{f(x)}{(1-ax)(1-bx)(1-cx)}$ has no pole at $\infty$. We use
the basic identity
 \begin{equation}\label{eqn:res}
    -\Res\left\{\frac{\,f(x)}{1-tx};\frac1t\right\}=\frac1t\,f\left(\frac1t\right)\,,
 \end{equation}
with $t=a\,,b\,,c$ respectively, for the calculation of residues
at a simple pole. %, supposing that $a$, $b$ and $c$ are such that
%$f$ does not have zeros at $\frac1a$, $\frac1b$ and $\frac1c$.
Substituting (\ref{eqn:res}) into (\ref{eqn:F}) we get
 $$
 F(a,b,c)=
    \frac{(a+u)^{g}(a+v)^{g}}{(a-b)(a-c)}+
    \frac{(b+u)^{g}(b+v)^{g}}{(b-a)(b-c)}+
    \frac{(c+u)^{g}(c+v)^{g}}{(c-a)(c-b)}\, .
 $$
Putting this into (\ref{Hodge-and-residue}), we get
 $$
    e(M(2,d)) = \frac{(1+u)^{g}(1+v)^{g}
                    (1+u^2v)^{g}(1+uv^2)^{g}
                    -(uv)^{g}(1+u)^{2g}(1+v)^{2g}}
                {(1-uv)(1-(uv)^2)}\,.
 $$
The computation only works for generic values of $u,v$ (because of
the restriction that we have imposed on $a,b,c$ above), but this
is enough to get the equality, since both terms are polynomials in
$u$ and $v$.
\end{dem}

\begin{flushleft}\obeylines
    Vicente Mu\~noz
    Departamento de Matem\'aticas
    Consejo Superior de Investigaciones Cient{\'\i}ficas
    Serrano 113 bis, 28006 Madrid, Spain
    \texttt{vicente.munoz@imaff.cfmac.csic.es}

    \

    Daniel Ortega
    Departamento de Matem\'aticas
    Facultad de Ciencias
    Universidad Aut\'onoma de Madrid
    28049 Madrid, Spain
    \texttt{daniel.ortega@uam.es}

    \

    Maria-Jes{\'u}s V{\'a}zquez-Gallo
    Departamento de Ingenier\'{\i}a Civil: Servicios Urbanos
    Unidad Docente: Matem\'aticas
    Escuela de Ingenier{\'\i}a de Obras P{\'u}blicas
    Universidad Polit{\'e}cnica de Madrid
    Alfonso XII 3 y 5, 28014 Madrid, Spain
    \texttt{mariajesus.vazquez@upm.es}

\end{flushleft}
\end{document}